\documentclass{article}

\usepackage{arxiv}

\usepackage[utf8]{inputenc} % allow utf-8 input
\usepackage[T1]{fontenc}    % use 8-bit T1 fonts
\usepackage{hyperref}       % hyperlinks
\usepackage{url}            % simple URL typesetting
\usepackage{booktabs}       % professional-quality tables
\usepackage{amsfonts}       % blackboard math symbols
\usepackage{nicefrac}       % compact symbols for 1/2, etc.
\usepackage{microtype}      % microtypography
\usepackage{graphicx}
\usepackage[square,sort,comma,numbers]{natbib}
\usepackage{doi}
\usepackage{float}
\usepackage{caption}
\usepackage{subcaption}
\usepackage{xcolor}
\usepackage{amsmath}
\usepackage{graphicx}
\usepackage{overpic}
\usepackage{xcolor}
\usepackage[dvipsnames]{xcolor}

\usepackage{listings}
\usepackage[section]{placeins}
\usepackage{float}

\title{A High-Order Spectral Element Solver for Steady-State Free Surface Flows}

\author{\href{https://orcid.org/0009-0009-4402-0998}{\includegraphics[scale=0.06]{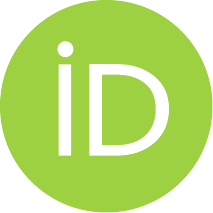}\hspace{1mm}Simone Minniti} \\
	Department of Applied Mathematics and Computer Science\\
	Technical University of Denmark\\
	Kongens Lyngby, 2800, Denmark \\
	\texttt{s232445@dtu.dk} \\
\And
 \href{https://orcid.org/0000-0001-6698-2623}{\includegraphics[scale=0.06]{Figures/orcid-eps-converted-to.pdf}\hspace{1mm}Jens Visbech} \\
	Department of Applied Mathematics and Computer Science  \\
	Technical University of Denmark\\
	Kongens Lyngby, 2800, Denmark \\
	\texttt{jvis@dtu.dk} \\
        \And 
    \href{https://orcid.org/0000-0001-6934-634X}{\includegraphics[scale=0.06]{Figures/orcid-eps-converted-to.pdf}\hspace{1mm}Claes Eskilsson} \\
    Maritime Department\\
    RISE - Research Institutes of Sweden \\
    SE-412 96 Gothenburg, Sweden \\
    \texttt{claes.eskilsson@ri.se}
     \\
        \And 
    \href{https://orcid.org/0000-0002-2497-7276}{\includegraphics[scale=0.06]{Figures/orcid-eps-converted-to.pdf} \hspace{1mm}Nicola Parolini} \\
	MOX Laboratory, Department of Mathematics \\
    Politecnico di Milano \\
    Milano, 20133, Italy\\
	\texttt{nicola.parolini@polimi.it}
\And
    	\href{https://orcid.org/0000-0001-8626-1575}{\includegraphics[scale=0.06]{Figures/orcid-eps-converted-to.pdf}\hspace{1mm}Allan Peter Engsig-Karup} \\
	Department of Applied Mathematics and Computer Science\\
	Technical University of Denmark\\
	Kongens Lyngby, 2800, Denmark \\
	\texttt{apek@dtu.dk} 
}

\renewcommand{\headeright}{Minniti et al. (2025)}
\renewcommand{\shorttitle}{An \textit{arXiv} Preprint}

\hypersetup{
pdftitle={},
pdfsubject={},
pdfauthor={},
pdfkeywords={},
}

\begin{document}
\maketitle

\begin{abstract}
	We present a spectral element solver for the steady incompressible Navier-Stokes equations subject to a free surface. Utilizing the kinematic behaviour of the free surface boundary, an iterative pseudo-time procedure is proposed to determine the a priori unknown free surface profile. The numerical model is implemented in the open-source finite element framework \texttt{Firedrake}, which enables the use of a high-order polynomial basis on unstructured meshes through weak formulations. Additionally, the curvature of the free surface and submerged bodies is incorporated through curvilinear elements obtained via transfinite linear blending, which conserves the high-order convergent properties of the overall scheme.
The model is applied to several benchmark cases in two spatial dimensions. Initially, it addresses fixed-domain problems, including the lid-driven cavity flow and flows around bodies such as a cylinder and a NACA airfoil. Subsequently, with the presence of a free surface, it is extended to determine the flow around a bathymetry bump and a submerged NACA airfoil. The results confirm the high-order accuracy of the model through convergence studies and demonstrate a substantial speed-up over low-order numerical schemes.
\end{abstract}

% keywords: (can be removed)

% Input your different sections here:
\section{Introduction}

Computational fluid dynamics (CFD) has become an indispensable tool in modern engineering, particularly for fluid-structure interaction modeling \cite{BonheureGazzolaSperone2019}, enabling the cost-effective analysis of complex fluid flows that are difficult, costly, or even impossible to study experimentally. Using numerical simulations to complement experiments can reduce the need for physical prototyping and accelerate research and design processes.

In this work, we consider the branch of CFD that deals with viscous free surface flows in the context of water waves. Our focus is on steady-state flow problems, whose analysis can be carried out at a reduced computational cost compared to fully transient simulations, as discussed in \cite{Chande2018Investigation, Fuzesi2024Importance}, where steady and transient approaches are compared and their respective advantages are thoroughly investigated. The study of steady-state solutions provides fundamental insight into the dominant flow structures and performance characteristics of engineered systems. In particular, we focus on the steady incompressible Navier–Stokes equations (INSE) which have been widely employed to model a broad range of viscous flow phenomena, including the wave patterns generated by submerged or partially submerged moving bodies. Both numerical techniques and experimental approaches have been developed to investigate this class of problems and to gain a deeper understanding of the underlying physical mechanisms.
For example, in \cite{nikitin2017adaptive}, an adaptive numerical method was presented for free surface flows past rigidly mounted obstacles, while in \cite{Wemmenhove2006numerical}, a compressible two-phase flow model was proposed for numerical simulation of hydrodynamic wave loading.
Moreover, mean drag and lift forces from steady-state CFD can be combined with spectral wave models to estimate fatigue life. For energy capture using wave energy devices \cite{Eskilsson2018WaveDragon}, steady-state CFD can help quantify mean flow rates, pressure differences, or velocity fields that are crucial to estimate energy conversion efficiency. From an experimental perspective, there are several works, such as \cite{Kim2023Experimental}, which presented a study on free surface waves generated by the motion of a circular cylinder-shaped submerged body in a single fluid layer, and \cite{Seong2022Experimental}, which investigated the influence of different velocities and accelerations on free surface waves produced by a moving body in a towing tank.

For the class of steady flows targeted in this work, we aim to address some of the key challenges associated with accurately and efficiently simulating them numerically.
It should be noted that open-source CFD software, e.g., \texttt{OpenFOAM} \cite{openfoam}, as well as commercial packages such as \texttt{STAR-CCM+} \cite{sratccm} and \texttt{ANSYS Fluent} \cite{ansysfluent}, mainly rely on low-order numerical schemes. Although these schemes are robust and widely applicable, achieving a desired level of accuracy often requires a substantially refined mesh, which can lead to considerably high computational costs for flows that require fine-scale resolution.
%Secondly, it is important to consider that in free surface flows, gravity plays a dominant role, which makes Froude similitude the natural basis for scaling between a prototype model and other real-world models \cite{Heller2011Scale}. However, enforcing the Froude similitude also constrains the Reynolds number, which characterizes the viscous flow regime around submerged bodies, such as NACA airfoils. As a result, while Froude scaling ensures that the wave-making and free surface adjustment are captured correctly, it often comes at the cost of mismatched Reynolds numbers, leading to differences in boundary-layer behavior, separation, and viscous drag. 

Moreover, in modeling free surface flows, the position of the free surface is a priori unknown and must be determined as part of the solution procedure. As the fluid flows past structural features or around a submerged body, the free surface continuously adapts in response to the evolving pressure and velocity fields in the fluid domain. To obtain a steady-state solution, an equilibrium configuration must be determined such that the velocity, pressure, and free surface profile are simultaneously balanced. Consequently, the INSE governing the fluid motion are strongly coupled with the dynamic and kinematic conditions at the free surface \cite{hino1988}.

%%% 2) Literature review: 
Various numerical approaches have been developed to address the free surface flow problem, each employing different techniques, sometimes tailored to slightly different formulations of the governing equations. Early hybrid approaches were proposed, such as the integral equation method based on Green’s theorem introduced in \cite{yeung1979}. Finite difference formulations were adopted in \cite{salvesen1976, coleman1986}, while finite element methods (FEM) were proposed in \cite{bai1978, wu1995, Lohner1999}. \cite{HIRTETAL1981, Ali2010Numerical} investigated the benefits of capturing the free surface profile using the volume-of-fluid (VOF) method, whereas the finite volume method (FVM) was employed to solve the steady-state flow in \cite{Hino1993, Karim2014Numerical}. Potential flow formulations for two-dimensional hydrofoils advancing at constant speed were also addressed in \cite{kennell1984, forbes1985}. Early work on the use of spectral element methods for viscous free surface problems was presented in \cite{HO1994207}. More recently, approaches such as the isogeometric analysis (IGA) method have been proposed in \cite{Akkerman2011Isogeometric}. A Reynolds-Averaged Navier–Stokes (RANS) solver, which captures both breaking and non-breaking waves generated by submerged hydrofoils, was presented in \cite{mascio2007}.
Experimental investigations have also played an essential role in understanding steady-state free surface flow problems. In particular, the measurements of free surface elevations and both breaking and non-breaking wave resistance for a NACA0012 hydrofoil at varying submergence depths, angles of attack, and velocities, presented in \cite{duncan1983}, remain a fundamental reference. 

% \cite{hino1988}.  \cite{kwag1991} \cite{bal1999}  SKIPPED 

%%% 3) Moving towards our work: 
In the setting of finite elements, which can use unstructured meshes for handling complex geometries encounted in real-world settings, we consider a spectral element method (SEM), a technique that combines the geometric flexibility of the FEM with the high accuracy of spectral methods. Due to its ability to deliver high-order approximations on complex geometries, SEM has become a well-established tool in the numerical solution of partial differential equations, particularly in fluid dynamics problems with original work due to \cite{patera1984spectral}, and since then it has been considered for free surface incompressible fluid flow problems in many works such as \cite{FischerEtAl1988, RobertsonSherwin1999, ENGSIGKARUP20161, EngsigKarupEtAl2019}. Popular textbooks describing SEM are \cite{karniadakis2005spectral, Quarteroni2017, kopriva2009implementing}, and a recent review of SEM is presented in \cite{Xu2018}. Moreover, open-source software platforms such as \texttt{FireDrake} \cite{FiredrakeUserManual} and \texttt{Nektar++} have been developed for the solution of partial differential equations, taking into account complex geometries using the spectral/$hp$ element method \cite{Cantwell2015Nektar, Moxey2020Nektar}.
% To accurately represent curved boundaries and avoid the introduction of dominant geometric errors within the computational domain, curvilinear elements are incorporated. This is achieved using a transfinite linear blending technique based on the methodologies outlined in \cite{Gordon1973, canuto2007fundamentals, visbech2022}. The approach enables precise geometric modelling, which is crucial for maintaining spectral convergence, which can otherwise degrade due to second-order geometric approximation errors near curved boundaries.
In the present work, to capture the free surface behaviour, an iterative approach similar to the one proposed in \cite{Lohner1999, VanBrummelen2001} is employed, in which the free surface boundary is adjusted at each pseudo-time step to converge to a steady-state solution. This dynamic process is crucial in obtaining an accurate representation of the steady-state flow and the wave pattern generated by the moving body. 

% STEADY-STATE CITATIONS
% An example of steady-state flow analysis can be found in \cite{Bensler1998Experimental}, where the flow generated by a spark-ignition engine intake port is investigated. The study focuses on comparing CFD and experimental results, with the aim of optimizing engine performance through steady-state CFD simulations.

% A density-based topology optimization framework for the design of three-dimensional heat sinks is proposed in \cite{Alexandersen2016Large}, wherein the steady INSE are coupled with the thermal convection-diffusion equation to accurately capture the flow and heat transfer behaviour. 

\subsection{Paper contributions}

The paper presents a SEM-based model of the INSE subjected to a free surface under steady-state conditions. The novelty of the work lies in the design of a velocity-pressure solver with inclusion of viscosity and in the use of high-order SEM for the development of a solver that can address steady-state free surface flow problems, features not present in previous studies. The new solver is implemented in \texttt{Firedrake}, an open-source finite element framework \cite{FiredrakeUserManual}, and will be demonstrated through a series of benchmark cases serving the purpose of verification and validation \cite{OBERKAMPF2002209}, including flows in fixed fluid domains (such as the lid-driven cavity and flow around objects in fixed domains), as well as flows over bathymetry changes and submerged bodies under free surface conditions. Considering these cases, we emphasize numerical validation against results available in the literature, assessment of numerical performance, and comparisons with conventional CFD simulations performed using \texttt{OpenFOAM}.
We demonstrate that, for fixed accuracy levels, a high-order SEM scheme achieves better computational efficiency than conventional low-order FEM schemes, with speed-ups of up to 210$\times$ in two space dimensions.
Hence, we propose high-order SEM as a basis for a new, powerful tool for numerical modeling of steady-state free surface flows, which is useful for engineers and researchers in this field.

\subsection{Outline of paper}

The remainder of this study is organized as follows. In Section \ref{sec:mathematicalproblem}, we introduce the governing equations for describing fluid flow in terms of the INSE with a free surface for steady-state problems. In Section \ref{sec:numericalapproach}, we detail the numerical scheme employed for the numerical discretization of the governing equations and the iterative approach used to reach the steady-state solution, where the elevation of the free surface is not prescribed a priori but must adapt consistently with the fluid flow. Details of the new implementation of weak forms in the open-source \texttt{Firedrake} framework are presented, followed by an explanation of the mesh generation procedure, the strategy for mesh updates, and the integration of curvilinear elements. In Section \ref{sec:Results}, we present a set of numerical experiments conducted for verification and validation purposes, starting from fixed-domain problems and concluding with more complex cases dealing with free surfaces. Enhancements were implemented in the Firedrake solver framework to include support for curvilinear elements based on isoparametric representations to effectively handle curved boundaries. The performance of low- and high-order methods is compared and discussed, with validation against established results from the literature. In Section \ref{sec:conclusion}, we conclude with a summary of our main results and suggestions for improvements.
\section{Mathematical Problem}
\label{sec:mathematicalproblem}

We consider a two-dimensional fluid domain, $\Omega \subset \mathbb{R}^2$, spanned by the spatial coordinates $\boldsymbol{x} = (x,y)$. The boundary of the domain consists of five distinct interfaces: inflow $\Gamma_{\text{in}}$, outflow $\Gamma_{\text{out}}$, free surface $\Gamma_\eta$, bed $\Gamma_{\text{bed}}$, and solid body $\Gamma_{\text{body}}$.
Since we restrict our analysis to two-dimensional configurations, we present the mathematical formulation of the problem in 2D, including the domain, the governing equations, and the boundary conditions. The extension to the three-dimensional setting is straightforward, as analogous forms of the equations and boundary conditions apply.

\begin{figure}[ht]
    \centering
    \includegraphics[width=0.9\textwidth]{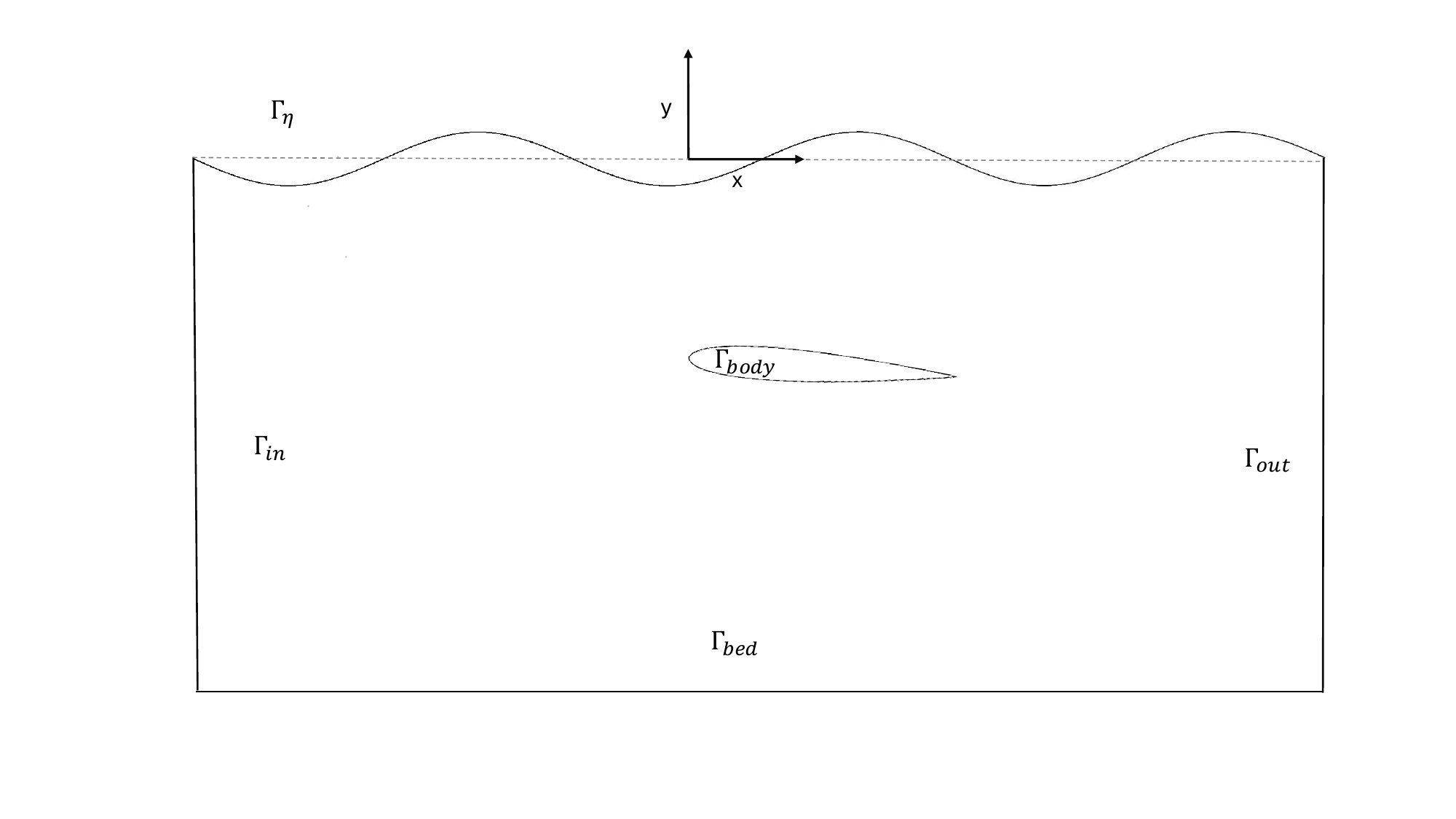}
    \caption{Computational domain and boundaries for the free surface flow over a submerged NACA0012 airfoil.}
    \label{fig:domain}
\end{figure}

The coordinate system is defined such that $y=0$ coincides with the undisturbed free surface level, while the choice of $x=0$ depends on the specific problem under consideration. For instance, in the case of a submerged body, $x=0$ is positioned at the leading edge of the body. A conceptual layout of the domain $\Omega$ and its boundaries is provided in Figure \ref{fig:domain}.

We introduce two key physical quantities: the free surface elevation, $\eta = \eta(x) : \Gamma_{\eta} \mapsto \mathbb{R}$, and the still-water depth, $h = h(x) : \Gamma_{\text{bed}} \mapsto \mathbb{R}$, measured from the seabed to $y=0$. 

\textbf{Remark.} For intermediate problems without a free surface such as the lid-driven cavity flow or channel flow problems, different domains are considered. Specifically, for the lid-driven cavity, a squared domain with a moving lid is chosen, while for the channel cases (NACA and cylinder), a rectangular domain with fixed boundaries is used. Details of domain sizes will be provided alongside the results.

\subsection{The steady incompressible Navier-Stokes equations}

Assuming the fluid to be homogeneous, steady, and incompressible, steady INSE can be employed to describe the flow within $\Omega$ in terms of the non-dimensional velocity field $\boldsymbol{u} = (u,v)  : \Omega \to \mathbb{R}^2$ and the non-dimensional pressure $p : \Omega \mapsto \mathbb{R}$, e.g., see \cite{Hino1989}. 
By targeting a steady-state solution, the flow variables are independent of time and thus depend solely on the spatial coordinates, i.e.,
\(
u = u(x,y), \ v = v(x,y), \ p = p(x,y).
\)

The governing steady INSE, incorporating the Newtonian stress tensor, are given by
\[
    \nabla \cdot \boldsymbol{u} = 0,\ \quad \text{and } \quad
    -\nabla \cdot \!\left( \frac{1}{\mathrm{Re}} \big( \nabla \boldsymbol{u} + (\nabla \boldsymbol{u})^{T} \big) \right)
    + (\boldsymbol{u} \cdot \nabla)\boldsymbol{u} 
    + \nabla p = \boldsymbol{f}.
\]

Here, the vector \(\boldsymbol{f}\) denotes additional body forces per unit volume, while the Reynolds number is defined as  
\(
    \mathrm{Re} = \frac{U L}{\nu},
\)  
where \(U\) and \(L\) are the characteristic velocity and length scale of the problem, respectively, and \(\nu\) is the kinematic viscosity of the fluid assumed to be constant. In a steady-state context, all flow quantities are independent of time, and consequently, no terms involving temporal derivatives appear. 

Moreover, assuming that gravity is the only external force acting on the free surface flow implies that the body force term reduces to the gravitational contribution alone, with no additional effects such as surface tension, density gradients, or other external fields. Consequently, its effect can be incorporated into the pressure field, following the formulation given in \cite{Lohner1999}. This is achieved by introducing the non-dimensional dynamic pressure $p_d : \Omega \mapsto \mathbb{R}$ as 
\[
p_d = p - \frac{-y}{\mathrm{Fr}^2},
\] 
where the term $\frac{-y}{\mathrm{Fr}^2}$ represents the non-dimensional hydrostatic pressure, assuming the vertical coordinate $y$ points upwards. The Froude number is defined as  
\(
    \mathrm{Fr} = \tfrac{U}{\sqrt{gL}},
\)  
with \(g\) denoting the gravitational acceleration.
This definition effectively removes the hydrostatic contribution from the total pressure, isolating the dynamic component.

By introducing the dynamic pressure \(p_d\) to account for the effect of gravity,  the governing equations in \(\Omega\) take the form

\begin{align}
\nabla \cdot \boldsymbol{u} = 0, \quad \text{and} \quad
-\nabla \cdot \!\left( \frac{1}{\mathrm{Re}} \big( \nabla \boldsymbol{u} + (\nabla \boldsymbol{u})^{T} \big) \right)
+ (\boldsymbol{u}\cdot\nabla)\boldsymbol{u}
+ \nabla p_d = 0. \label{eq:NSE}
\end{align}

We impose specific boundary conditions on the domain boundaries, following the approach outlined in \cite{Hino1989,Hino1993}. A uniform horizontal velocity \(\boldsymbol{u} = (U, 0)\) is prescribed on $\Gamma_{\text{in}}$, valid under the assumption of a fixed water level at the inflow boundary. A no-slip condition, \( \boldsymbol{u} = \boldsymbol0 \), is imposed on \( \Gamma_{\text{body}} \), while an impermeability condition is imposed on \( \Gamma_{\text{bed}} \), requiring the normal component of the velocity to vanish, i.e., \( \boldsymbol{u} \cdot \boldsymbol{n} = 0 \). Additionally, considering that for an incompressible Newtonian fluid the Cauchy stress tensor is \( \boldsymbol{\sigma} = -p \mathbf{I} + \frac{1}{\mathrm{Re}} \left( \nabla \boldsymbol{u} + (\nabla \boldsymbol{u})^T \right) \), it follows that the stress vector acting on a boundary with normal \( \boldsymbol{n} \) is \( \boldsymbol{T} = \boldsymbol{\sigma} \cdot \boldsymbol{n} = \left( -p \mathbf{I} + \frac{1}{\mathrm{Re}} \left( \nabla \boldsymbol{u} + (\nabla \boldsymbol{u})^T \right) \right) \cdot \boldsymbol{n} \). 

Since the bed does not impose any friction, we require that the tangential component of the stress vanishes, which leads to the condition
\(
\left(
\frac{1}{\mathrm{Re}}
\big( \nabla\boldsymbol{u} + (\nabla\boldsymbol{u})^{T} \big)\boldsymbol{n}
- p_d\,\boldsymbol{n}
\right)\!\cdot\boldsymbol{t} = \boldsymbol{0} 
\  \text{on } \Gamma_{\text{bed}} .
\)
Finally, on \( \Gamma_{\text{out}} \) we consider the natural Neumann condition
\(
\frac{1}{\mathrm{Re}}
\big( \nabla\boldsymbol{u} + (\nabla\boldsymbol{u})^{T} \big)\boldsymbol{n}
- p_d\,\boldsymbol{n}
= \boldsymbol{0} .
\) This condition imposes a stress-free condition on the outlet boundary and naturally arises in the weak formulation of the problem; further details are provided in Section \ref{sec:weak}. Lastly, the kinematic and dynamic boundary conditions that hold on \( \Gamma_\eta \) are introduced in the following section.

\subsection{Kinematic and dynamic free surface boundary conditions}

Due to the presence of a free surface, INSE must be coupled with the kinematic and dynamic free surface boundary conditions at $\Gamma_{\eta}$, which govern the evolution of the free surface elevation.

On the free surface $\Gamma_\eta$, a constant pressure level can be assumed, i.e., \(
p = 0 \ \text{on } \Gamma_\eta.
\)
Consequently the dynamic pressure on the free surface, ${p}_{d}$, is given by
\begin{equation}
{p}_{d} = \frac{\eta}{\mathrm{Fr}^2},  \quad \boldsymbol{x} \in \Gamma_{\eta},
\label{eq:dynamicCondition}
\end{equation}
which will be referred to as the dynamic free surface boundary condition.

It is important to note that this boundary condition is not imposed strongly, but rather enters naturally through the weak formulation of the governing equations. At the free surface, continuity of stresses must hold. Assuming that surface tension effects are negligible, the dynamic boundary condition reduces to a Neumann-type condition in which the fluid stress tensor balances the prescribed free surface pressure, corresponding to a stress-free interface. In particular, this leads to

\begin{equation}
\frac{1}{\mathrm{Re}} \big( \nabla \boldsymbol{u} + (\nabla \boldsymbol{u})^{T} \big) \boldsymbol{n}- p_d \boldsymbol{n} = -{p}_{d} \boldsymbol{n}, \implies 
\frac{1}{\mathrm{Re}}
\big( \nabla \boldsymbol{u} + (\nabla \boldsymbol{u})^{T} \big)\boldsymbol{n}
= \boldsymbol{0} \quad \boldsymbol{x} \in \Gamma_{\eta}.
\label{eq:Neumann2}
\end{equation}

Moreover, the impermeability of the free surface is enforced through the kinematic condition, its derivation is detailed in \cite{Minniti2025}. Considering that $ y = \eta(x)$ is the surface level, the equation results in a pure advection equation as

\begin{equation}
\frac{\partial \eta}{\partial t} + u \frac{\partial \eta}{\partial x}  = v, \quad \boldsymbol{x} \in \Gamma_{\eta}.
\label{eq:kinematicCondition}
\end{equation}

Equation \eqref{eq:kinematicCondition}, referred to as the one-dimensional kinematic free surface condition (i.e. defined on a manifold of one spatial dimension less than the fluid domain), plays a crucial role in the steady-state pseudo-time-stepping procedure described in Section \ref{sec:numericalapproach}. It is important to emphasize that the time derivative is not null, as this condition is actively employed within the iterative pseudo-time process.
\section{Numerical Discretization}
\label{sec:numericalapproach}

In this section, we first present the weak formulation of the problem, followed by the Galerkin approximation and its implementation in \texttt{Firedrake}. After establishing how the INSE and the free surface conditions are discretized and solved, we then introduce the iterative strategy employed to compute the steady-state solution of the free surface problem.

\subsection{Weak formulations} 
\label{sec:weak}

Considering the equations and boundary conditions presented in Section \ref{sec:mathematicalproblem}, we now derive the weak formulation of the Navier-Stokes problem \eqref{eq:NSE}. 
The objective is to determine the velocity field 
\( \boldsymbol{u} \in V = [H^1(\Omega)]^2 \) and the dynamic pressure 
\( p_d \in Q = L^2(\Omega) \) that satisfy the system in a variational sense. 

To formulate the problem variationally, we introduce a suitable test space \( V_0 \subset [H^1(\Omega)]^2 \) which denotes the subspace of vector fields that vanish on \( \Gamma_{\text{in}} \cup \Gamma_{\text{body}} \) 
and have a zero normal component on \( \Gamma_{\text{bed}} \). 
The weak formulation is obtained by multiplying the governing equations by the test functions \( \boldsymbol{v} \in V_0 \) and \( q \in Q \), then integrating over the domain, and applying Green's identities. In this way, the Neumann boundary conditions are naturally incorporated. Moreover, consider that 
\(
\mathbf{D}(\boldsymbol{u}) = \frac{1}{2}\big( \nabla \boldsymbol{u} + (\nabla \boldsymbol{u})^{T} \big)
\)
denotes the deviatoric tensor, the resulting weak problem reads as follows.

Find \( \boldsymbol{u} \in [H^1(\Omega)]^2 \) with 
\( \boldsymbol{u} = \boldsymbol{g} \) on \( \Gamma_{\text{in}} \cup \Gamma_{\text{body}} \) 
and \( \boldsymbol{u} \cdot \boldsymbol{n} = 0 \) on \( \Gamma_{\text{bed}} \), and 
\( p_d \in Q = L^2(\Omega) \), such that 
\begin{equation}
\begin{aligned}
   & a(\boldsymbol{u},\boldsymbol{v}) + c(\boldsymbol{u}, \boldsymbol{u},\boldsymbol{v}) 
+ b(\boldsymbol{v},p_d) = F(\boldsymbol{v}) \quad \forall \boldsymbol{v} \in V_0, \\
& b(\boldsymbol{u},q) = 0 \quad \forall q \in Q,
\end{aligned}
\label{eq:weakproblem}
\end{equation}

with the linear, bilinear and trilinear forms defined as

\begin{align*}
    a(\boldsymbol{u},\boldsymbol{v}) 
   & = \frac{2}{\mathrm{Re}} \int_{\Omega} \mathbf{D}(\boldsymbol{u}) : \mathbf{D}(\boldsymbol{v}) \, d\Omega, 
    & b(\boldsymbol{v},q) &= - \int_{\Omega} q \, \nabla \cdot \boldsymbol{v} \, d\Omega, \notag \\
    c(\boldsymbol{u}, \boldsymbol{w}, \boldsymbol{v}) &= \int_{\Omega} (\boldsymbol{u} \cdot \nabla) \boldsymbol{w} \cdot \boldsymbol{v} \, d\Omega, 
    & F(\boldsymbol{v}) &= - \int_{\Gamma_{\eta}} {p}_{d} \,\boldsymbol{n} \cdot \boldsymbol{v} \, d\Gamma,
    \label{eq:forms}
\end{align*}

considering the identity
\(
\mathbf{D}(\boldsymbol{u}) : \nabla \boldsymbol{v}
= \mathbf{D}(\boldsymbol{u}) : \mathbf{D}(\boldsymbol{v}).
\)

It is worth noting that, while the dynamic surface condition \eqref{eq:dynamicCondition} enters the problem as a Neumann-type boundary condition, instead the kinematic condition \eqref{eq:kinematicCondition} plays a central role in the iterative approach. Consequently, deriving its variational formulation and selecting an appropriate numerical pseudo-time stepping scheme are essential for the stability and accuracy of the computation. As explained in detail in Section \ref{sec:iterative}, we employ an iterative procedure based on a pseudo-time variable, denoted by $\tau$. Moreover, we consider a finite difference scheme to approximate the time derivative in  \eqref{eq:kinematicCondition} before deriving the weak scheme \footnote{This is a Rothe's method where temporal discretization is done before spatial discretization of differential equations.},
\(
\frac{\partial \eta}{\partial t} \approx \frac{\eta_{\tau+1}-\eta_\tau}{\Delta \tau},
\)
where $\eta_{\tau+1}$ and $\eta_\tau$ denote the elevations of the free surface in consecutive pseudo-time steps, with $\Delta \tau$ chosen according to the CFL condition \cite{Quarteroni2014}.

At the inflow boundary, a Dirichlet condition fixes the height of the water column \(\eta(x_{\text{in}})=0 \ \forall \tau \) and it's imposed to ensure constant flux. Although this assumption may initially seem restrictive because it does not permit adjustment of the surface at the boundary, we ensure that it would have a negligible impact on the flow dynamics above and downstream of the object, provided that the inflow is positioned sufficiently far from the object.

Multiplying the modified equation \eqref{eq:kinematicCondition} by a test function \(\phi \in \Phi = H^1_0 (\Gamma_\eta)\) and integrating over the free surface spatial domain $\Gamma_\eta \subset \mathbb{R}$ we obtain the weak formulation of the problem: find \(\eta_{\tau+1} \in H^1  (\Gamma_\eta)\) with \(\eta_{\tau+1}(x_{\text{in}})=0\) such that
\begin{equation}
\int_{\Gamma_\eta} \left( \frac{\eta_{\tau+1}-\eta_\tau}{\Delta \tau} 
+ u_\tau \frac{\partial \eta_{\tau+1}}{\partial x}  - v_\tau \right)\phi \, dx = 0, 
\quad \forall \phi \in \Phi.
\label{eq:freesurfaceWF}
\end{equation}

It is important to note that for a given pseudo-time step $\tau$, the only unknown in \eqref{eq:freesurfaceWF} will be $\eta_{\tau+1}$. Indeed, when solving this equation, the velocity field will be fixed and known. The complete numerical procedure will be presented in more detail in Section \ref{sec:iterative}. It is noteworthy that $\frac{\partial \eta_{\tau+1}}{\partial x}$ is used instead of $\frac{\partial \eta_{\tau}}{\partial x}$ to adopt an implicit scheme for improved numerical temporal stability for pseudo-time stepping.

\textbf{Remark.} Following \cite{Lohner1999}, an artificial damping term is incorporated into \eqref{eq:freesurfaceWF} in the vicinity of the outflow boundary. This addition prevents wave reflections and ensures the numerical stability of the simulation.  

Specifically, the free surface kinematic condition restricted to one space dimension is modified to include a damping term (cf. \cite{Engsig2013Designing}) as:

\begin{equation}
\frac{\partial \eta}{\partial t} + u \frac{\partial \eta}{\partial x} - v + \gamma\eta = 0,
\end{equation}
where \(\gamma = \gamma(x)\) represents the artificial damping term, defined as:

\begin{equation}
\gamma(x) =
\begin{cases} 
% A \left( \frac{x - x_{d1}}{x_{in} - x_{d1}} \right)^2, & \text{if } x_{in} \leq x \leq x_{d1},  \\
% A \left( \frac{x - x_{d2}}{x_{out} - x_{d2}} \right)^2, & \text{if } x_{d2} \leq x \leq x_{out}, \\
A \left( \frac{x - x_{d}}{x_{out} - x_{d}} \right)^2, & \text{if } x_{d} \leq x \leq x_{out}, \\
0, & \text{otherwise}.
\end{cases}
\label{eq:damping}
\end{equation}

Here, \(A\) is an arbitrary constant that determines the strength of the damping, typically set as $A=C/d\tau$ where $C=\mathcal{O}(1)$ and $d\tau$ is the time-stepping size \cite{Engsig2013Designing}. The parameter \(x_{d}\), which defines the extent of the damping region, is given by
\(
 x_{d} = x_{out} - 2\pi Fr^2.
\)

% Thus, the damping region is positioned at a distance corresponding approximately to one wavelength, as estimated by linear theory, from the outflow boundary \cite{Hino1993}. The quadratic formulation of the damping term presented in Equation \eqref{eq:damping} ensures a gradual increase in dissipation as the inflow and outflow boundaries are approached, effectively mitigating wave reflections.

\subsubsection{Galerkin approximation}

To obtain a numerical solution for the weak formulation, a Galerkin approximation based on finite-dimensional subspaces must be introduced. 
Let us consider the polygonal domain $\Omega$. We can introduce a computational grid $\mathcal{T}_h$ and two finite-dimensional subspaces \( V_h \subset V\) and \( Q_h \subset Q \). Consider $V_h^0 = V_h \cap V_0$, where we recall that $Q$, $V_0$ and $V$ have been defined in Section \ref{sec:weak}. 
The Galerkin method then consists of restricting the infinite-dimensional problem to these finite-dimensional spaces. 

The Galerkin approximation of problem \eqref{eq:weakproblem} reads: find \( (\boldsymbol{u}_h, p_{d_h}) \in {V}_h \times Q_h \), with \( \boldsymbol{u}_h = \boldsymbol{g}_h \text{ on } \Gamma_{\text{in}} \cup \Gamma_{\text{body}} \text{ and }\boldsymbol{u}_h \cdot \boldsymbol{n} = 0 \text{ on } \Gamma_{\text{bed}} \), such that for all \( (\boldsymbol{v}_h, q_h) \in V_h^0 \times Q_h \):
\begin{equation}
\begin{aligned}
    a(\boldsymbol{u}_h, \boldsymbol{v}_h) + c(\boldsymbol{u}_h, \boldsymbol{u}_h, \boldsymbol{v}_h) + b(\boldsymbol{v}_h, p_{d_h}) &= F(\boldsymbol{v}_h),  \\
    b(\boldsymbol{u}_h, q_h) &= 0, 
    \label{eq:galerkin}
\end{aligned}
\end{equation}
where $\boldsymbol{g}_h$ is an approximation of the Dirichlet data $\boldsymbol{g}$ in the space $V_h( \Gamma_{\text{in}} \cup \Gamma_{\text{body}})$. Note that it is important to consider $V_h$ and $Q_h$ such that they satisfy the discrete inf-sup condition. For sufficiently small data and for subspaces which satisfy this condition, it is possible to demonstrate that the solution of the Galerkin approximation \eqref{eq:galerkin} converges to the solution of the original continuous problem \eqref{eq:weakproblem} \cite{Quarteroni2014}. 

% The well-posedness of the discrete problem \eqref{eq:galerkin1}–\eqref{eq:galerkin2} can be established analogously to the continuous case, provided that the inf-sup condition \eqref{eq:discrete-inf-sup} holds.

% There exist several possible choices for the discrete velocity and pressure spaces,  \( V_h \) and \( Q_h \), which must be selected so as to satisfy the discrete inf-sup condition. A straightforward approach would be to employ equal-order polynomial spaces for both the velocity and pressure fields; however, such a choice generally violates this condition. 

Consider the definition of the polynomial space \( \mathbb{P}_N \):

\[
\mathbb{P}_N(\Omega) 
= \left\{\, p : \Omega \to \mathbb{R} \;\middle|\; 
p(\boldsymbol{x}) = \sum_{|\boldsymbol{\alpha}| \leq N} 
c_{\boldsymbol{\alpha}}
\, \boldsymbol{x}^{\boldsymbol{\alpha}}, \;
c_{\boldsymbol{\alpha}} \in \mathbb{R} \,\right\},
\]

where \(\boldsymbol{\alpha} = (\alpha_1, \dots, \alpha_d)\) is a multi-index with \(|\boldsymbol{\alpha}| = \sum_{i=1}^d \alpha_i\) and \(\boldsymbol{x}^{\boldsymbol{\alpha}} = \prod_{i=1}^d x_i^{\alpha_i}\).

To satisfy the discrete inf-sup condition and ensure compatibility at higher order, we employ the mixed Taylor-Hood finite element pair
\(\mathbb{P}_k / \mathbb{P}_{k-1}\) for \(k \geq 2\), defining the discrete velocity and pressure spaces \(V_h\) and \(Q_h\), respectively.

\subsection{Implementation of the spectral element solver in Firedrake}
\label{sec:firedrake}

The resulting formulation is directly implemented in \texttt{Firedrake} \cite{FiredrakeUserManual}. 
\texttt{Firedrake} integrates with PETSc, enabling efficient parallel computations \cite{petsc-web-page, balay2024petsc}, which is essential for large-scale problems requiring effective utilization of computational resources. Additionally, it leverages the Unified Form Language (UFL) \cite{Alnaes2014UFL} to define variational forms, thereby simplifying the implementation of complex weak formulations.

Previous studies demonstrate both the advantages of spectral element methods and the flexibility of the computational framework. In \cite{Jacobs2015Firedrake}, a numerical model of shallow water flows was presented using \texttt{Firedrake}; \cite{Bokhove2016Variational} introduced an efficient \texttt{Firedrake} solver for the Benney–Luke-type equations; more recently, \cite{visbech2025fnpfsemparallelspectralelement} employed the framework to perform parallel simulations of water waves and their interactions with offshore structures using a fully nonlinear potential flow (FNPF) model. Furthermore, the Thetis project, an unstructured grid coastal ocean model, was built using this framework \cite{Karna2018Thetis}.  

Extending the weak formulation introduced in Section \ref{sec:weak} we define the residual \(R\) as

\begin{equation}
\begin{aligned}
R = \frac{2}{\mathrm{Re}} \int_{\Omega} \mathbf{D}(\boldsymbol{u}) : \mathbf{D}(\boldsymbol{v}) \, d\Omega + \int_{\Omega} (\boldsymbol{u} \cdot \nabla) \boldsymbol{u} \cdot \boldsymbol{v} \, d\Omega \ - \int_{\Omega} p_d \, \text{div} \boldsymbol{v} \, d\Omega + \int_{\Omega} q \, \text{div} \boldsymbol{u} \, d\Omega \ +\int_{\Gamma_{\eta}} p_{d} \boldsymbol{n} \cdot \boldsymbol{v} \ d\Gamma.
 \label{eq:weakNS2nocrossflow}
\end{aligned}
\end{equation}
% PENALTY
% + K \int_{\Gamma_{\text{bed}}} (\boldsymbol{u} \cdot \boldsymbol{n})(\boldsymbol{v} \cdot \boldsymbol{n}) \, d\Gamma.

\texttt{Firedrake} allows the definition of function spaces, residuals, and boundary conditions in a precise and concise manner. To obtain a nodal high-order spectral element method, we employ Lagrange basis functions to define the finite-dimensional spaces of the desired order on the constructed mesh. We denote by $P$-$Q$ the polynomial degree of the pair of finite element spaces for velocity and pressure, respectively. For instance, when using $P$-$Q$ = $2$-$1$ indicates that the velocity field is approximated in the quadratic space $\mathbb{P}_{2}$, while the pressure field is approximated in the linear space $\mathbb{P}_{1}$. Dirichlet boundary conditions are imposed by specifying the velocity values directly on the relevant boundaries identified by specific tags.

In this work, triangular unstructured meshes are used. For the triangular elements, \texttt{Firedrake} employs a set of non-equidistant nodes that share properties with standard Gauss-Legendre-Lobatto (GLL) (p) nodes, coinciding with them along the edges. Their construction, detailed in \cite{Tobin2020}, is based on a blending strategy between higher-dimensional entities, combining the features of the GLL(\(k\)) nodes for \(k = 1, \dots, p\).

The nonlinear problem defined by the residual \eqref{eq:weakNS2nocrossflow} is solved using the \texttt{solve} function, which relies on the \textit{PETSc} framework \cite{petsc-web-page} and employs a Newton-based iterative method. At each iteration, the Jacobian of the residual form is assembled via automatic differentiation of the residual with respect to the solution variable. In general, the zero field can be used as the initial guess for the Newton iterations. However, for free surface flow simulations, the solution from the previous pseudo-time step is used as the initial guess, as the flow changes incrementally with the temporal updates.
The Jacobian system within each Newton step is solved using PETSc's default Krylov subspace method, GMRES, with a strict relative tolerance of $10^{-8}$ and an absolute tolerance of $10^{-10}$; in Section \ref{sec:Results} we verify that the proposed solver is working. Preconditioning is provided by PETSc’s built-in incomplete LU factorization; however, in ongoing work, it may be subject to improvement in performance and scalability by employing defect correction with multi-grid techniques as in related works \cite{EngsigKarup2014, EngsigKarupEtAl2021, Melander2024pmultigridacceleratednodalspectral} and using new subspace acceleration techniques for improved initial guesses \cite{guido2024subspaceaccelerationsequencelinear,sorensen2025subspaceaccelerationefficientnonlinear}. 

\subsection{Mesh generation}

Meshes are generated using \texttt{gmsh} \cite{gmsh}, allowing refinement in regions with expected high velocity and pressure gradients, and coarser discretization elsewhere. Depending on domain size, Reynolds number, orders involved, and accuracy requirements, multiple meshes are created and refined through uniform subdivisions to enhance resolution while preserving mesh topology. While uniform structured meshes are used for simpler problems, such as the lid-driven cavity, unstructured meshes are employed for more complex cases. 
For example, for the NACA airfoil problem, unstructured meshes with finer elements near the airfoil are used to capture the steep velocity gradient. An example of such a mesh is shown in Figure \ref{fig:mesh_naca_channel}. Observe that a dynamic adaptation algorithm is considered to preserve the quality of the mesh while the free surface evolves; see Section \ref{sec:meshupdate} for more details.

\subsection{Iterative approach}
\label{sec:iterative}

\begin{figure}[ht]
    \centering
\includegraphics[width=0.9\textwidth]{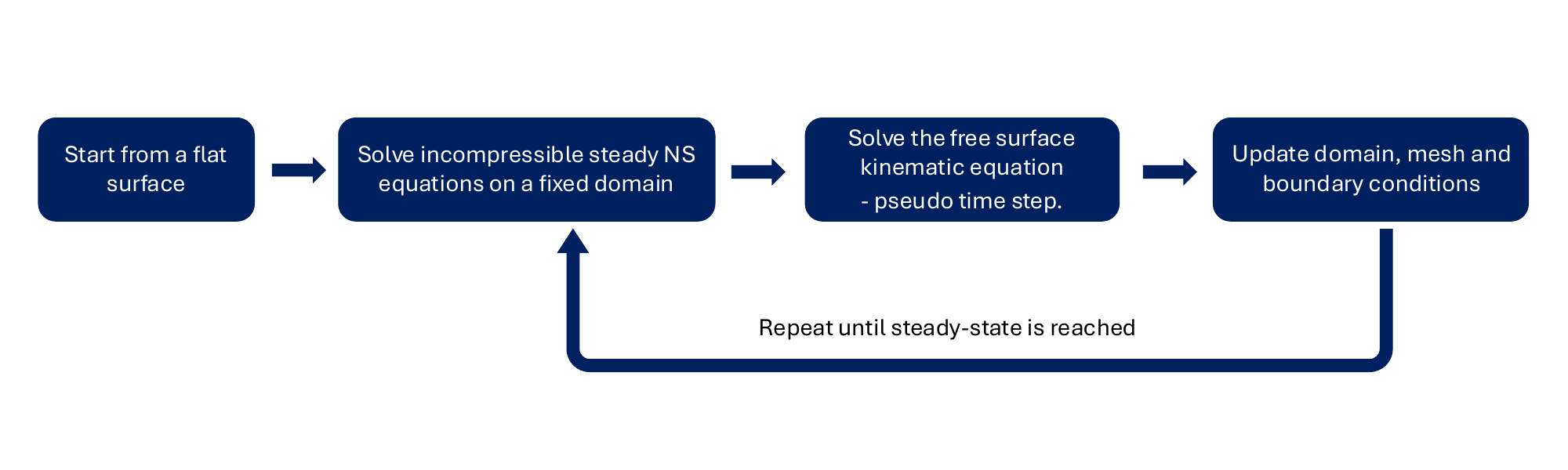}
    \caption{The iterative numerical procedure to reach steady-state.}
    \label{fig:num}
\end{figure}

The numerical procedure aims to compute the steady-state solution of a free surface flow using an interface-tracking strategy inspired by previous works \cite{Lohner1999,Hino1993}. The iterative procedure based on pseudo-timestepping is initialized with a flat free surface profile. Then, the INSE are solved at each step in a fixed computational domain to obtain the velocity and pressure fields. The velocity components are then extracted along the current free surface profile and used to solve the kinematic boundary equation \eqref{eq:freesurfaceWF}, where the only unknown is $\eta_{\tau+1}$. 
In this way the procedure provides the updated free surface profile, defining an updated physical fluid domain. The mesh and the corresponding boundary conditions are then consequently adapted. The surface nodes are displaced vertically according to the new elevation, and the deformation is propagated to the interior nodes using a mesh smoothing technique that ensures element quality is preserved (fx. avoid significant stretching of elements). More details are provided in Section \ref{sec:meshupdate}. Once the mesh has been updated, the dynamic boundary free surface condition is redefined \eqref{eq:dynamicCondition}, the residual \eqref{eq:weakNS2nocrossflow} is consequently updated, and the INSE are solved again in the new domain. The procedure is iterated until convergence is reached, as illustrated in Figure \ref{fig:num}.

Convergence is evaluated by monitoring the relative displacement in the free surface elevation between successive pseudo-time iterations, defined as
\begin{equation}
D = \frac{\| \eta_{\tau+1} - \eta_{\tau} \|_2}{\| \eta_{\tau} \|_2},
\label{eq:residual}
\end{equation}
% where \(\eta_{\tau}\) and \(\eta_{\tau+1}\) denote the free surface profiles at pseudo-time levels \(\tau\) and \(\tau+1\), respectively.

Specifically, the criterion is considered satisfied once $D$ falls below a prescribed tolerance, typically set to $10^{-6}$, indicating that a steady-state configuration has been reached.

\subsection{Mesh update}
\label{sec:meshupdate}

The entire mesh is updated during each iteration to prevent the formation of excessively large or distorted elements during the displacement of the free surface nodes. Specifically, once the free surface elevation at the next pseudo-time step, \(\eta_{\tau+1}\), is computed, the mesh is consistently adapted to the new configuration.
This transformation aims to smoothly propagate the displacement from the free surface to the underlying nodes, preserving overall mesh quality. 

\begin{figure}[ht]
    \centering
    \begin{subfigure}{0.48\textwidth}
        \centering
    \includegraphics{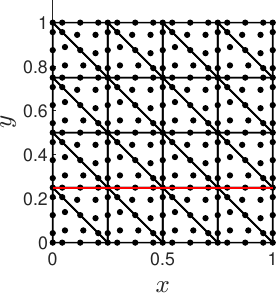}
        \caption{Initial mesh, $y_0 = 0.25$.}
    \end{subfigure}
    \begin{subfigure}{0.48\textwidth}
        \centering
        \includegraphics{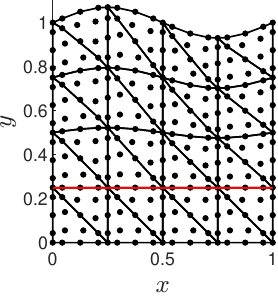}
        \caption{Mesh updated with $\eta(x) = C\sin(2\pi x)$.}
    \end{subfigure}
    
    \caption{Illustration of the mesh updating strategy using sinusoidal free surface displacements with amplitude \(C = 0.07\). The lower vertical threshold \(y_0\) is highlighted in red.}

    \label{fig:mesh_update}
\end{figure}

First, the \( y \)-coordinates of the free surface nodes are updated using the newly computed elevation values \( \eta_{\tau+1} \).  
Next, the horizontal domain is partitioned into bins, with the number of bins matching the number of free surface points. For each bin, the mesh nodes whose \( x \)-coordinates fall within its range are identified. Among these, only the nodes located above the immersed body or the bathymetry bump are considered, with the parameter \(y_0\) defining a lower vertical threshold.
For each selected node, the \( y \)-coordinate is updated at each iteration according to the difference between the new and the previous free surface elevations in the corresponding bin. The new vertical position is computed as
\begin{equation}
y_{\text{new}} = y_{\text{old}} + \min \left(1, \frac{y_{\text{old}} - y_0}{\eta_{\tau+1}^* - y_0} \right)(\eta^*_{\tau+1} - \eta_\tau^*),
\end{equation}

where \( \eta^* \) denotes the free surface elevation associated with the current bin.  
This scaling ensures that nodes closer to the free surface are displaced more significantly, while deeper nodes experience minor adjustments, leading to a smooth mesh deformation. All \( x \)-coordinates are maintained fixed throughout the entire procedure. Figure \ref{fig:mesh_update} illustrates how the mesh is updated in response to variations in the free surface curvature.  
A polynomial order of $P=4$ is considered, so each element contains 15 degrees of freedom (DOFs). As a result, curvilinear elements are generated due to the distortion of the mesh.  
This procedure allows for a more accurate representation of the surface profile, even when using relatively coarse meshes.

% FULL 4 FIGURES

% \begin{figure}[H]
%     \centering
%     \begin{subfigure}{0.48\textwidth}
%         \centering
%     \includegraphics{Figures/mesh01.pdf}
%         \caption{Initial mesh, $y_0 = 0.25$.}
%     \end{subfigure}
%     \begin{subfigure}{0.48\textwidth}
%         \centering
%     \includegraphics{Figures/mesh02.pdf}
%         \caption{Mesh updated with $\eta(x) = \tfrac{1}{3}C\sin(2\pi x)$.}
%     \end{subfigure}

%     \vspace{0.4cm}

%     \begin{subfigure}{0.48\textwidth}
%         \centering
%     \includegraphics{Figures/mesh03.pdf}
%         \caption{Mesh updated with $\eta(x) = \tfrac{2}{3}C\sin(2\pi x)$.}
%     \end{subfigure}
%     \begin{subfigure}{0.48\textwidth}
%         \centering
%     \includegraphics{Figures/mesh04.pdf}
%         \caption{Mesh updated with $\eta(x) = C\sin(2\pi x)$.}
%     \end{subfigure}
% \end{figure}
\subsection{Curvilinear elements}
\label{sec:curvilinearElements}

As anticipated, one of the key advantages of the SEM is its geometric flexibility, which enables an accurate representation of complex domains. However, when curvilinear boundaries are present, a second-order geometric approximation error may dominate if the geometry is represented using straight-edged elements. As discussed in \cite{Quarteroni2014}, increasing the polynomial degree of the basis functions can significantly reduce the approximation error (achieving spectral convergence). However, this advantage is lost if the geometry is not accurately captured. For this reason, it is essential to adopt curvilinear elements to provide a more faithful approximation of the boundaries, e.g., see \cite{canuto2007spectral,canuto2007fundamentals, Moxey201541}.

% First, it is important to note that only the elements located on or adjacent to curved boundaries require curvilinear treatment. The initial step thus consists of accurately identifying these elements and determining which edges must be curved. 

We present the methodology for constructing curvilinear isoparametric elements in two space dimensions using transfinite interpolation with linear blending, following the approaches described in \cite{Gordon1973, kopriva2009implementing}. This transformation proceeds by first expressing each degree of freedom in barycentric coordinates relative to the physical triangle, which are then used to locate the corresponding position in the canonical reference triangle. Once this consistent mapping between the physical and reference configurations is established, the blended Gordon-Hall mapping \cite{Gordon1973} is applied using the parametric definition of the curved boundary.

\begin{figure}[H]
    \centering\includegraphics{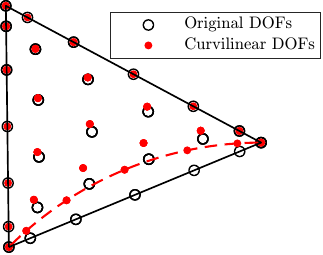}
    \caption{Curvilinear vs. straight-sided adapted-GLL nodes position on a $P=6$ triangular element. }\label{fig:curved_vs_straight}
\end{figure}

In Figure \ref{fig:curved_vs_straight}, the positions of the original DOFs are compared with those of the curvilinear element for a single triangular element with adapted-GLL nodes. 
It can be observed that the DOFs along edges that approximate curved boundaries are projected onto the curve, while those at vertices and on straight edges remain fixed. The internal DOFs are adjusted through a linear mapping based on their relative distance from the boundary.

Within Firedrake, this procedure is implemented by iterating over all elements and selecting those that contain an edge defining the boundary of the curvilinear object. For each selected element, the coordinates of its vertices are identified, ensuring a counter-clockwise ordering. As described above, all DOFs are then mapped first to the reference configuration and subsequently to the new physical curvilinear element. 
It is worth noting that, in \texttt{Firedrake}, the relocated degrees of freedom are reassigned directly to the mesh, after which the function spaces and element Jacobians are updated to account for the deformation.
\section{Numerical Experiments}
\label{sec:Results}

This section presents a series of benchmarks and numerical results for the steady incompressible Navier-Stokes solver. The aim is twofold: first, to validate the solver on established fixed-domain test cases and to assess the advantages of high-order discretizations over low-order ones; second, to present results for steady free-surface flow configurations, verifying the accuracy and robustness of the proposed solution strategy. To improve readability, the results are organized into two main subsections: fixed domain flows and free surface channel flows.

\subsection{Fixed domain flows}

We begin by validating the solver for confined, fixed-domain configurations. Specifically, we consider the lid-driven cavity flow and the flow around a NACA airfoil in a channel to assess accuracy and convergence properties. Once the base solver is validated, we analyze the impact of curvilinear elements in the simulation of flow around a circular cylinder. In this context, we also highlight the relative advantages of high-order spectral element schemes compared to low-order discretizations.

\subsubsection{Base solver validation: lid-driven cavity flow}

\begin{figure}[H]
\centering
% --- Re = 1000 ---
\begin{subfigure}{.5\textwidth}
  \centering
  \includegraphics{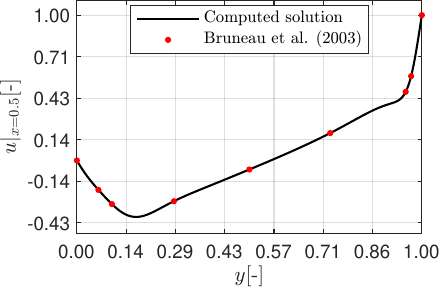}
  \caption{Horizontal velocity at $\text{Re} = 1{,}000$.}
  \label{fig:LiddrivenRe=1000_u}
\end{subfigure}%
\begin{subfigure}{.5\textwidth}
  \centering
  \includegraphics{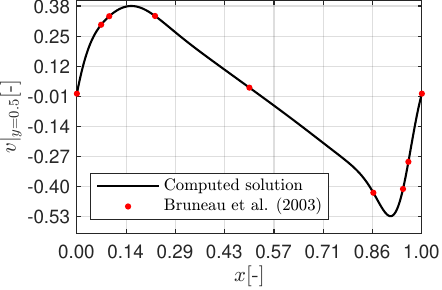}
  \caption{Vertical velocity at $\text{Re} = 1{,}000$.}
  \label{fig:LiddrivenRe=1000_v}
\end{subfigure}

% --- Re = 25000 ---
\begin{subfigure}{.5\textwidth}
  \centering
  \includegraphics{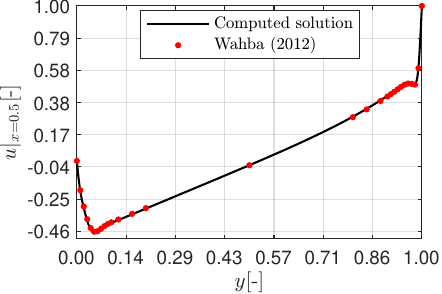}
  \caption{Horizontal velocity at $\text{Re} = 10{,}000$.}
  \label{fig:LiddrivenRe=25000_u}
\end{subfigure}%
\begin{subfigure}{.5\textwidth}
  \centering
  \includegraphics{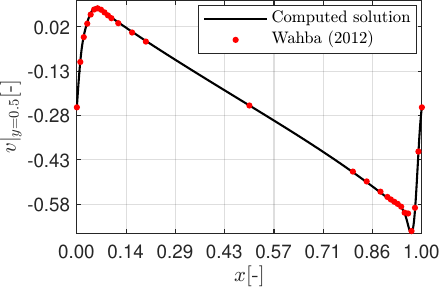}
  \caption{Vertical velocity at $\text{Re} = 10{,}000$.}
  \label{fig:LiddrivenRe=25000_v}
\end{subfigure}

\caption{Comparison between SEM results and the ones presented in \cite{WAHBA201285, Bruneau2006lid} for $\text{Re} = 1{,}000$ and $\text{Re} = 10{,}000$.}
\label{fig:comparisonLidDriven}
\end{figure}

We first solve the lid-driven cavity problem in the classical unit square domain $\Omega = (0,1)\times(0,1)$, with the top boundary moving at a constant velocity $\mathbf{u} = [1,0]$ [-]. We present validation cases in which the vertical velocity along the mid-horizontal axis and the horizontal velocity along the mid-vertical axis are computed for different Reynolds numbers, using a uniform structured triangular grid of resolution $128 \times 128$ with $P$-$Q$ = $3$-$2$. Figures \ref{fig:LiddrivenRe=1000_u} and \ref{fig:LiddrivenRe=1000_v} compare our results at $\text{Re}=1{,}000$ with those of \cite{Bruneau2006lid}, obtained using a third-order scheme. Figures \ref{fig:LiddrivenRe=25000_u} and \ref{fig:LiddrivenRe=25000_v} show the comparison at $\text{Re}=10{,}000$ with the reference steady solutions from \cite{WAHBA201285}, calculated with a compact fourth-order central difference scheme. The results align with established benchmarks and demonstrate that the proposed solver remains stable in the simulation of the steady INSE, even at high Reynolds numbers. It should be noted, however, that although a numerical steady solution is obtained for large Reynolds, it does not necessarily correspond to the physical flow, which is expected to become unsteady for $\mathrm{Re} > 1000$, with time-dependent perturbations developing.

\textbf{Remark.}
To obtain solutions at high Reynolds numbers, a continuation technique is employed. The problem is first solved starting from a null solution at a low Reynolds number (e.g., $\mathrm{Re}=1$). The resulting fields then serve as initial conditions for progressively higher values. Iterating this process enables us to achieve a stable solution even at very high Reynolds numbers. This approach is applied to all the other problems discussed in the following sections.

\subsubsection{Base solver validation: channel flow past a submerged NACA0012 airfoil}

\begin{figure}[H]
    \centering
    \begin{subfigure}[b]{0.65\textwidth}
        \centering
        \includegraphics[width=\textwidth]{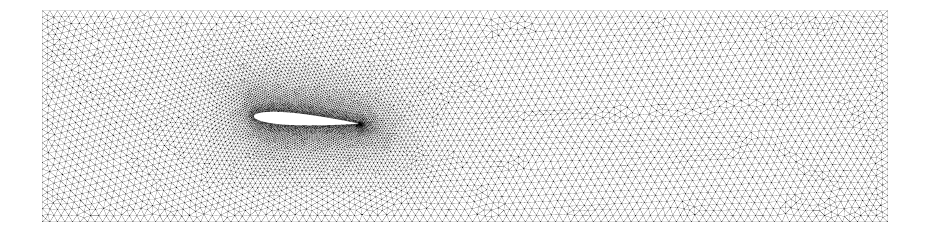}
        \caption{Unstructured mesh for a NACA0012 airfoil placed in a channel domain.}
        \label{fig:mesh_naca_channel}
    \end{subfigure}
    \vspace{0.4cm}
    \begin{subfigure}[b]{0.49\textwidth}
        \centering
        \includegraphics[width=\textwidth]{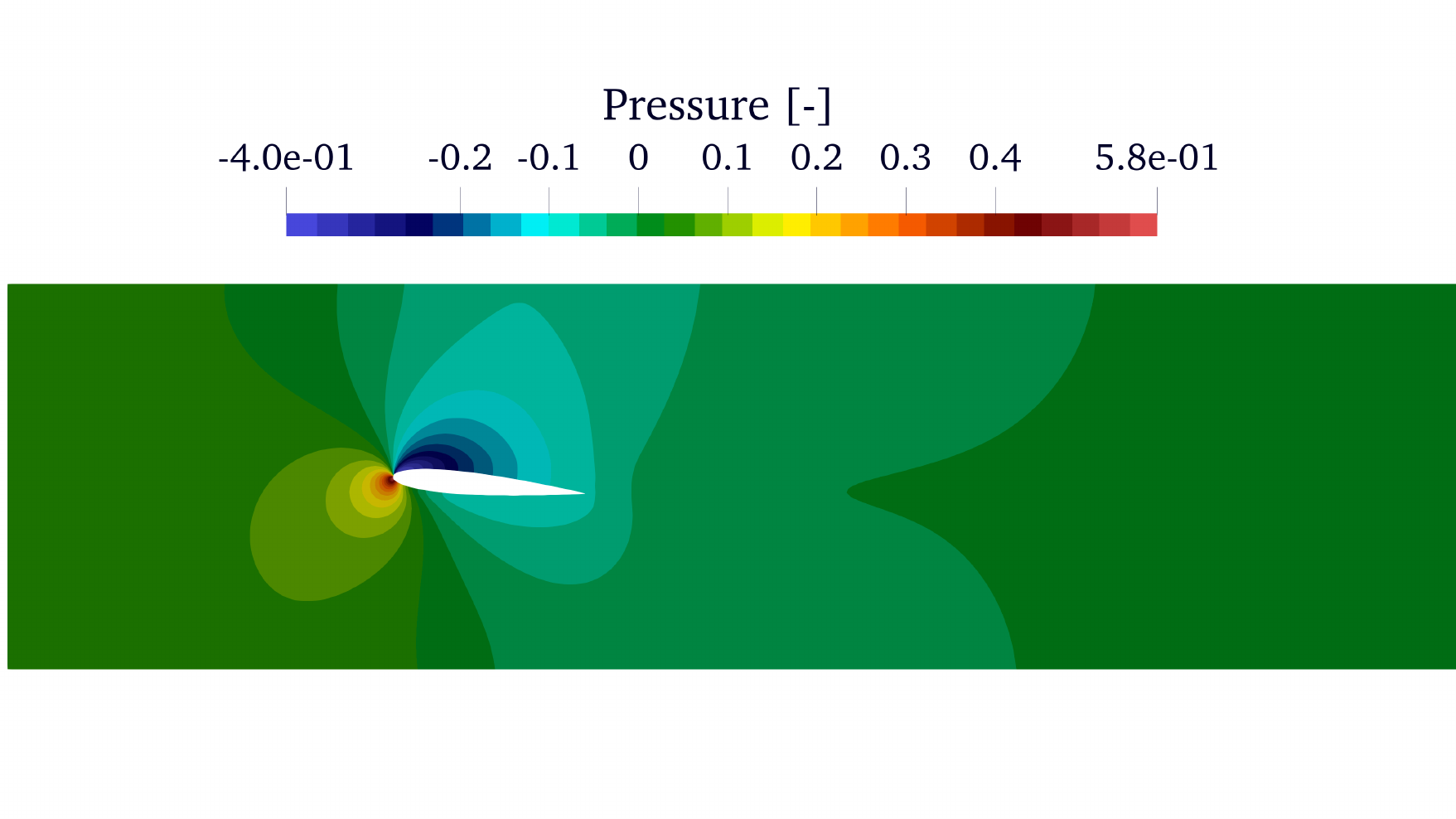}
        \caption{Dynamic pressure, \( \text{Re} = 1000 \).}
        \label{fig:pressure_naca_channel}
    \end{subfigure}
    \hfill
    \begin{subfigure}[b]{0.49\textwidth}
        \centering
        \includegraphics[width=\textwidth]{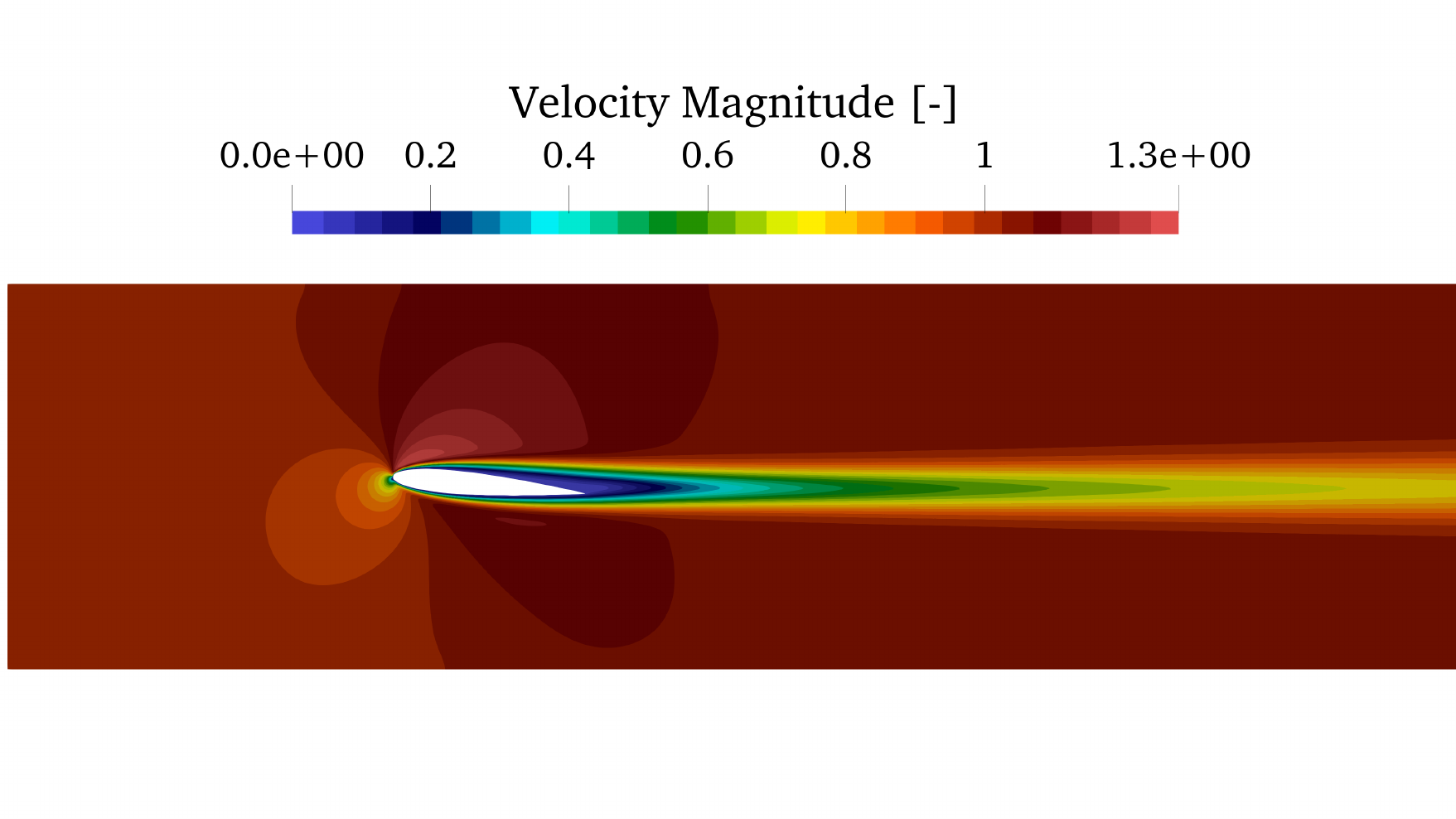}
        \caption{Velocity, \( \text{Re} = 1000 \).}
        \label{fig:velocity_naca_channel}
    \end{subfigure}
    \caption{Flow around a NACA0012 airfoil placed in a channel at \( \text{Re} = 1000 \). For visualization purposes, the results are computed on a reduced computational domain \( \Omega = [-2, 6] \times [-1, 1] \).}
    \label{fig:naca_channel_combined}
\end{figure}

Next, we consider the case of a fully submerged NACA0012 airfoil placed in a channel. In the computational domain, defined as \(\Omega = [-12,\,24] \times [-8,\,8]\), the airfoil, with chord length equal to 1 [-], is positioned such that its leading edge lies at the origin \((0,0)\). 
We denote by $\alpha$ its angle of attack. The domain configuration is chosen consistently with \cite{DIILIO2018200} to enable a fair comparison. An unstructured triangular mesh is generated and locally refined in the vicinity of the airfoil to adequately resolve the boundary layer.  
The resulting discretization comprises 10,714 nodes and 21,428 elements over the entire domain, with 200 edges used to represent the body. For $\text{Re}=1{,}000$, different angles of attack are investigated. The results show a symmetric flow field when the airfoil is not inclined, while larger recirculation vortices appear as the angle of attack $\alpha$ increases, as expected. In Figure \ref{fig:naca_channel_combined}, the mesh used for a smaller computational domain is shown (top), along with the corresponding velocity and pressure field solutions (bottom) obtained for \( \text{Re} = 1000 \), $P$-$Q$ = $3$-$2$, and an angle of attack of \( 5^{\circ} \).

\begin{figure}[H]
\centering
\begin{subfigure}{0.48\textwidth}
  \centering
  \includegraphics{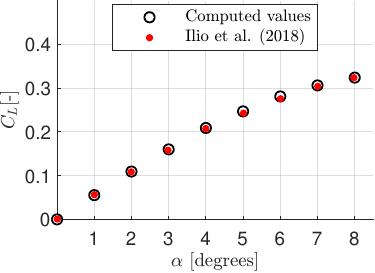}
  \caption{Lift coefficient.}
\end{subfigure}%
\begin{subfigure}{.48\textwidth}
  \centering
  \includegraphics{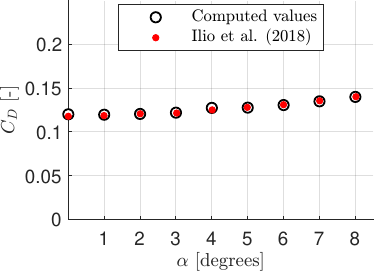}
  \caption{Drag coefficient.}
\end{subfigure}
\caption{Comparison between the computed drag and lift coefficients and the reference results from \cite{DIILIO2018200} at different angles of attack $\alpha$ for \(\text{Re} = 1{,}000\).}
\label{fig:CdClcomp}
\end{figure}

The non-dimensional lift and drag coefficients are computed choosing $P$-$Q$ = $3$-$2$ for different angles of attack. Results are compared with the benchmark results reported in \cite{DIILIO2018200}. The comparison between the computed values and the reference data is illustrated in Figure \ref{fig:CdClcomp}, demonstrating good agreement. 

\subsubsection{Curvilinear verification: channel flow past a submerged cylinder} 

In this section, we first demonstrate that the curvilinear elements constructed using the procedure described in Section \ref{sec:curvilinearElements} effectively mitigate the dominant geometrical error arising from curved geometries.
We then evaluate their impact on relevant fluid dynamic quantities, namely the drag and lift coefficients \cite{DIILIO2018200}, for the flow around a two-dimensional cylinder.

The computational domain is defined as a rectangle $\Omega = [-2,4] \times [-2,2]$ containing a circular cylinder of diameter $D=1$ [-] centred at the origin $(0,0)$. To evaluate the impact of curvilinear elements, a sequence of unstructured triangular meshes is generated and several combinations of polynomial degrees $P$-$Q$ are considered. The number of nodes, elements, and boundary elements is reported in Table \ref{tab:mesh_stats}.

\begin{table}[ht]
    \centering
    \caption{Mesh characteristics for different resolutions and polynomial orders.}
    \label{tab:mesh_stats}
    \begin{tabular}{crrrrr}
        \toprule
        Mesh & \# nodes & \# elements & \# boundary elements \\
        \midrule
        Mesh 0   & 72        & 144        & 8     \\
        Mesh 1   & 260       & 520        & 16    \\
        Mesh 2   & 984       & 1,968      & 32   \\
        Mesh 3   & 3,824     & 7,648      & 64   \\
        Mesh 4   & 11,829    & 23,658     & 128 \\
        Mesh 5   & 46,511    & 93,022     & 256  \\
        Mesh 6   & 185,148   & 370,296    & 512   \\
        Mesh 7   & 721,100   & 1,442,200  & 1,024\\
        \bottomrule
    \end{tabular}
\end{table}

\begin{figure}[ht]
    \centering
    \begin{subfigure}[b]{0.32\textwidth}
        \centering
        \includegraphics[width=\textwidth]{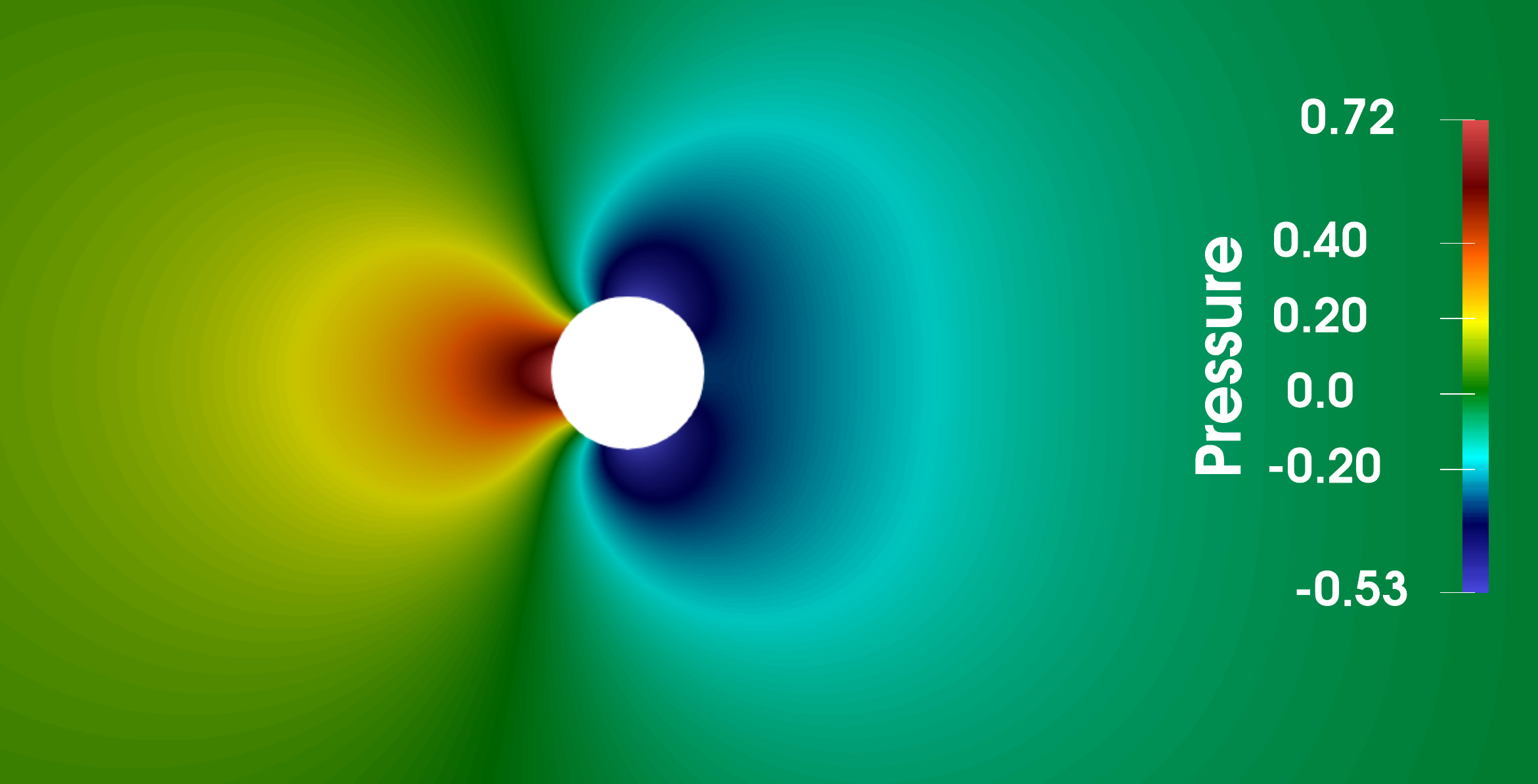}
        \caption{Dynamic pressure.}
        \label{fig:cylinder_pressure_re26}
    \end{subfigure}
    \hfill
    \begin{subfigure}[b]{0.32\textwidth}
        \centering
        \includegraphics[width=\textwidth]{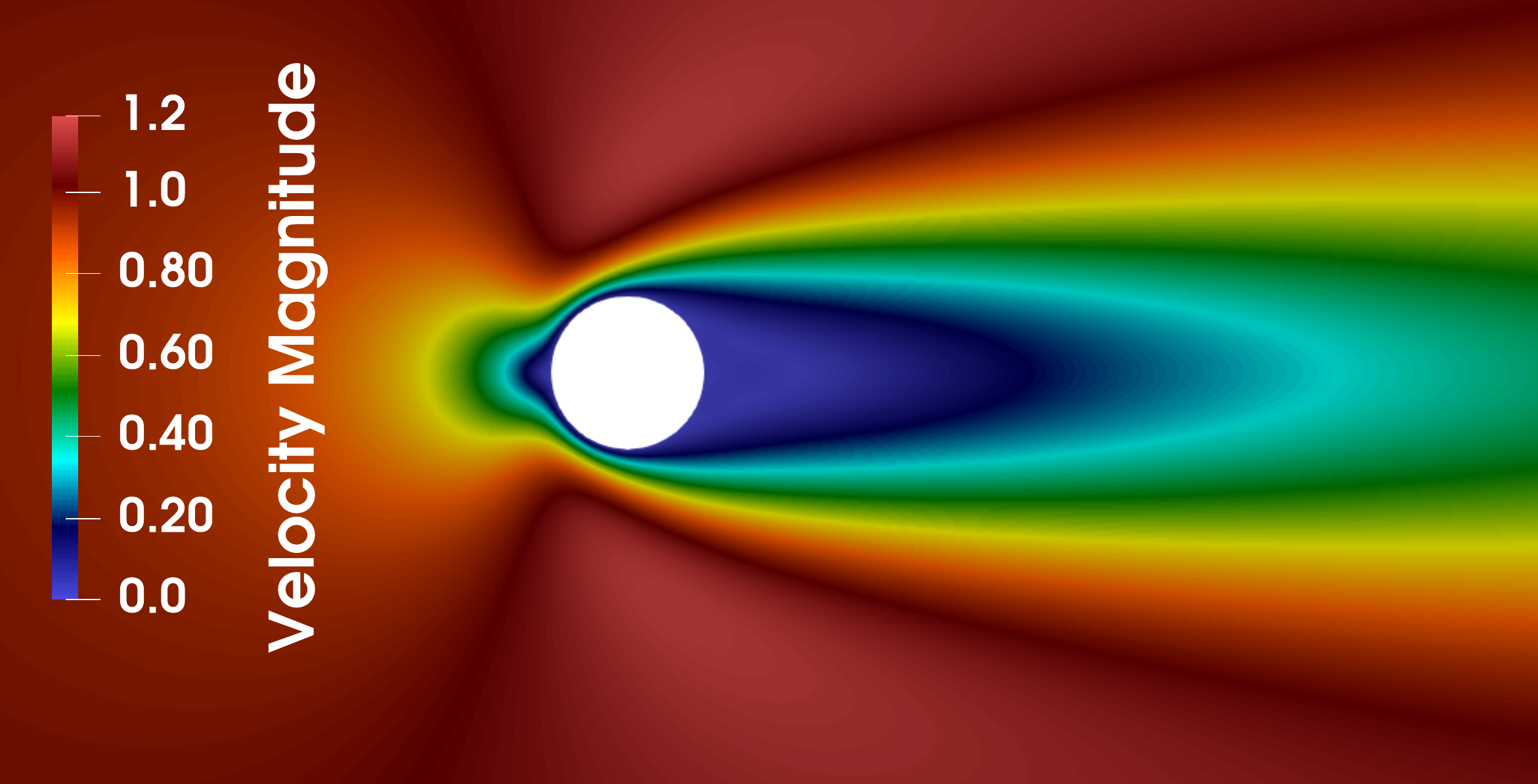}
        \caption{Velocity magnitude.}
        \label{fig:cylinder_streamlines_re26}
    \end{subfigure}
    \hfill
    \begin{subfigure}[b]{0.32\textwidth}
        \centering
        \includegraphics[width=\textwidth]{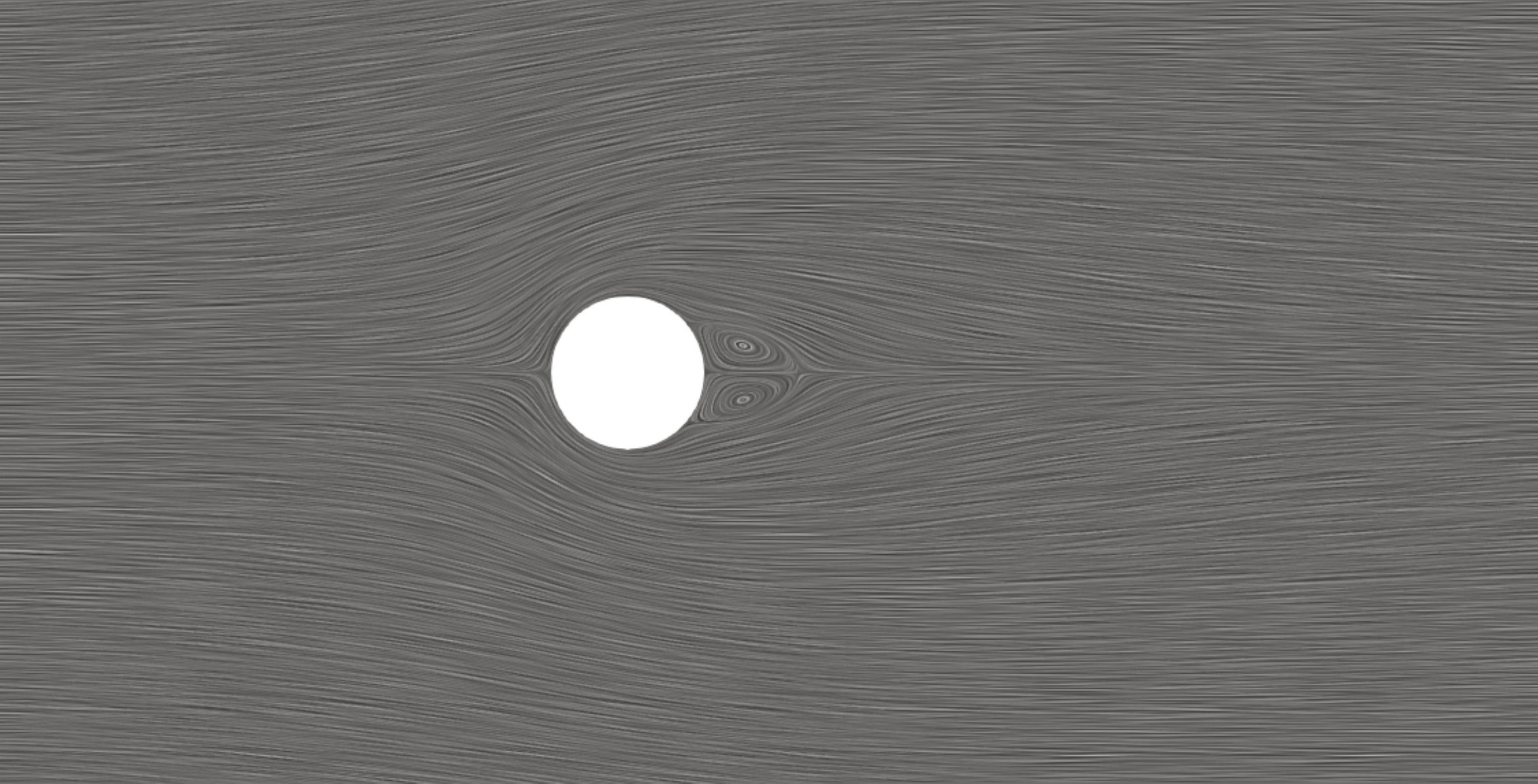}
        \caption{Streamlines.}
        \label{fig:cylinder_velocity_re26}
    \end{subfigure}
    \caption{Flow fields around a cylinder for $\mathrm{Re}=15$ using Mesh 3 with $P$-$Q=4$-$3$. Shown are dynamic pressure (a),  velocity magnitude (b), and streamlines (c).}
    \label{fig:cylinderFlow}
\end{figure}

Numerical results, such as those presented in Figure \ref{fig:cylinderFlow}, are qualitatively compared with experimental observations reported by \cite{VanDyke1982album} for $\mathrm{Re}=15$. Both experimental and numerical solutions exhibit a steady and symmetric separation bubble, which is characteristic of laminar flow for $4 < \mathrm{Re} < 40$.

\paragraph{Geometrical approximations.}

\begin{figure}[ht]
    \centering
    % First row
    \begin{subfigure}[b]{0.48\linewidth}
        \centering
        \includegraphics{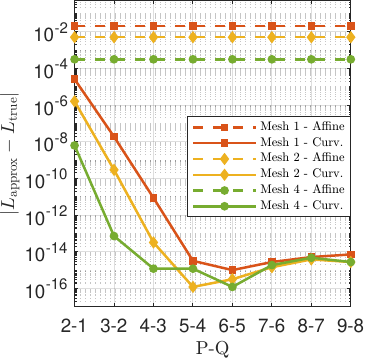}
        \caption{Convergence of the circular perimeter.}
        \label{fig:circle_perimeter_convergence}
    \end{subfigure}
    \hfill
    \begin{subfigure}[b]{0.48\linewidth}
        \centering
        \includegraphics{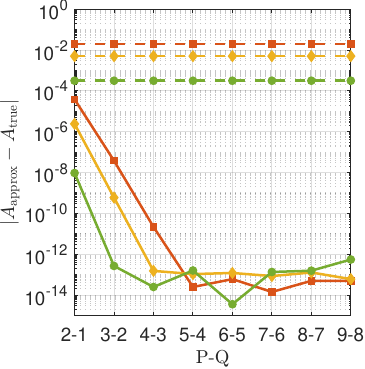}
        \caption{Convergence of the circular domain area.}
        \label{fig:circle_area_convergence}
    \end{subfigure}

    \vspace{0.5cm} % space between rows

    % Second row
    \begin{subfigure}[b]{0.48\linewidth}
        \centering
        \includegraphics{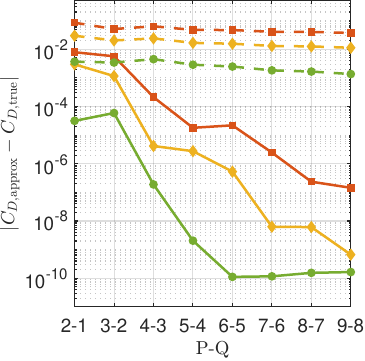}
        \caption{Drag coefficient convergence.}
        \label{fig:drag_convergence}
    \end{subfigure}
    \hfill
    \begin{subfigure}[b]{0.48\linewidth}
        \centering
        \includegraphics{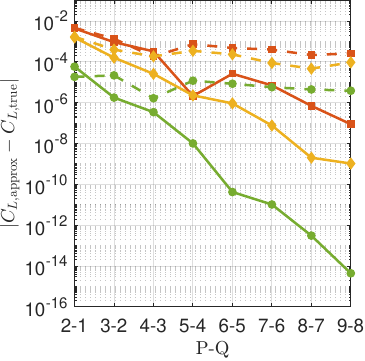}
        \caption{Lift coefficient convergence.}
        \label{fig:lift_convergence}
    \end{subfigure}

\caption{Comparison between affine and curvilinear elements for increasing polynomial orders and mesh refinements.}

    \label{fig:circle_geom_drag_lift}
\end{figure}

The geometric accuracy of the curved boundary is evaluated by computing the arc perimeter and the circular domain area through integration, ensuring that the Jacobian is consistently reassembled rather than relying on the moved coordinates. Figures \ref{fig:circle_perimeter_convergence} and \ref{fig:circle_area_convergence} show that the absolute errors in length and area decrease with increasing polynomial degrees \(P\)-$Q$ and mesh refinement. The affine elements exhibit constant errors (dashed lines), while the curvilinear elements achieve spectral convergence, eventually limited by the precision of the machine around \(10^{-15}\). These findings confirm that curvilinear elements provide a more accurate representation of the curvature of the geometry.  

\paragraph{Lift and drag approximations.}

For the analysis of drag and lift coefficients, we consider a uniform flow around a cylinder at a Reynolds number of $\mathrm{Re} = 15$. A visualization of the pressure and velocity fields is provided in Figure \ref{fig:cylinderFlow}. Note that the smallest mesh, Mesh 1 with order $2$-$1$, involves $260$ DOFs, whereas the largest one, Mesh 4 with order $9$-$8$, contains $942{,}021$ DOFs.

The error behavior for both straight-sided and curvilinear elements is compared in Figures \ref{fig:drag_convergence} and \ref{fig:lift_convergence}. For the lift coefficient, the error is measured as the absolute deviation from zero, since theory predicts it to vanish for this configuration. For the drag coefficient, it's defined as the absolute difference with respect to a reference value computed on a very fine high-order mesh, namely Mesh 6 with polynomial order $5$-$4$ including $4{,}610{,}780$ DOFs. As observed for the geometric quantities, both mesh refinement and higher polynomial degrees reduce the error, with spectral convergence apparent for curvilinear elements. Notably, without the transfinite blending procedure, increasing the polynomial degree has only a minor effect: indeed, although the number of DOFs increases, the convergence rate is ultimately limited by the geometric approximation when straight-sided affine elements are used.

\paragraph{Performance assessment: high-order vs low-order numerical schemes.}

Performance assessments were conducted on the Technical University of Denmark (DTU) high-performance computing (HPC) servers to investigate the computational benefits of high-order methods compared to low-order schemes. For each combination of mesh (see Table \ref{tab:mesh_stats} for mesh characteristics) and polynomial degree, multiple runs were performed, and the average CPU time was recorded. As expected, the computational time increases with both the polynomial degree and mesh refinement. Consider that for this case, the largest mesh involved is Mesh 7 with order $3$-$2$, containing $6{,}479{,}148$ DOFs.

To evaluate and compare the performance of different configurations in a fair way, we adopt the following approach: for four fixed error tolerances and for each mesh, the minimum polynomial degree required to meet the given tolerance is identified. A given tolerance is considered satisfied when the error in the lift coefficient falls below that threshold. \textbf{Remark.} Analogous results can be obtained when considering the drag coefficient.

\begin{figure}[ht]
    \centering
    \includegraphics{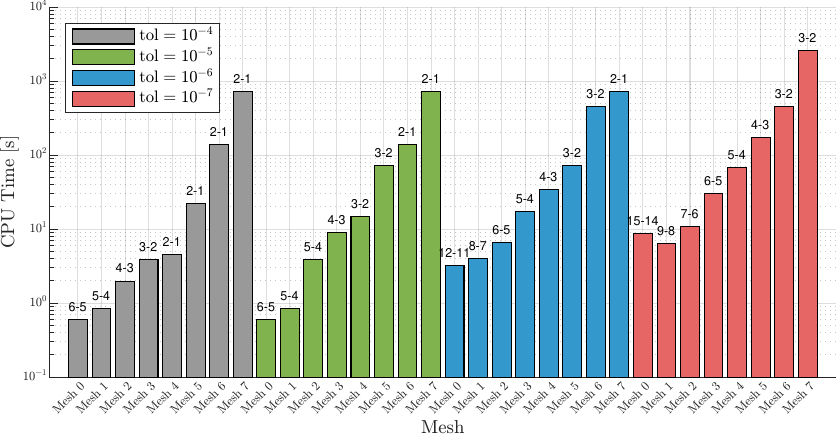}
    \caption{CPU time required to achieve fixed error tolerances in $C_L$. Labels indicate the lowest polynomial degrees $P$-$Q$ required to pass the tolerance.}
    \label{fig:ho_performance}
\end{figure}

Naturally, from Figure \ref{fig:circle_geom_drag_lift} it is obvious that coarser meshes demand higher polynomial degrees, whereas finer meshes may be able to achieve the same accuracy with lower-degree polynomials. The corresponding CPU times for the selected configurations are reported in Figure \ref{fig:ho_performance}. Above each bar in the figure, the corresponding polynomial degree combination is indicated. 

It is important to note that, as expected, each tolerance can be achieved either through a low-order scheme on a fine mesh or via a higher-order scheme on a coarser mesh. However, the key insight from the plot is that high-order methods on coarse meshes consistently outperform low-order methods on fine meshes in terms of computational time. Indeed, it is generally more efficient to use a coarse mesh with a high-order discretization than to rely on extensive mesh refinement to achieve the same accuracy. Observe that the CPU time on the vertical axis is displayed using a logarithmic scale. For instance, if a tolerance of $10^{-5}$ is considered, using a low-order scheme ($P$-$Q$ = $2$-$1$) requires six levels of mesh refinement (Mesh 6 is the coarsest that can achieve this tolerance with a degree of $2$-$1$) and approximately 137 seconds of computational time. In contrast, the same tolerance can be achieved in just $0.6$ seconds using Mesh 0 and a high-order scheme ($P$-$Q$ = $6$-$5$). 

\begin{figure}[ht]
    \centering
    \includegraphics{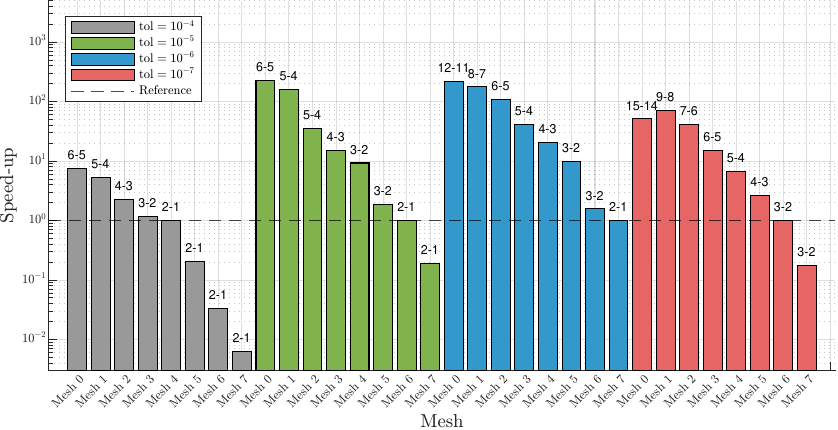}
    \caption{Speed-up for each tolerance level across different meshes, computed with respect to the fastest low-order case.}
    \label{fig:hosu}
\end{figure}

Figure \ref{fig:hosu} shows the speed-up, defined as $\text{speed-up} = t_{\text{reference}}/t_{\text{CPU}}$, where $t_{\text{reference}}$ is the computational time of the fastest lowest-order mesh for the corresponding tolerance.
For tolerances \(10^{-4}\), \(10^{-5}\), and \(10^{-6}\), the CPU time of the $2$-$1$ coarsest mesh at each tolerance is taken as the reference time. At tolerance \(10^{-7}\) Mesh 6 with polynomial order \(3\text{-}2\) achieves the required accuracy, with its CPU time serving as the reference time.
 Depending on the tolerance, the speed-up ranges from $7\times$ to $210\times$, clearly demonstrating the advantage of employing a high-order SEM over a low-order scheme at comparable accuracy.

Figure \ref{fig:ho_performance} and \ref{fig:hosu} clearly demonstrate the potential benefits of using high-order solvers for the class of problems considered in this work.

\subsection{Free-surface channel flow}

We now address steady free-surface flow problems by employing the numerical iterative procedure introduced in Section~\ref{sec:iterative}. Two benchmark cases are considered: steady flow over a bathymetry bump and flow past a fully submerged NACA0012 airfoil. We first present and analyze the results, then demonstrate the use of curvilinear elements and compare the resulting steady-state solutions with reference results computed using \texttt{OpenFOAM} in similar domain settings.

\subsubsection{Free surface flow validation: submerged bump}

In the context of free surface flows, we first consider a uniform flow over a submerged bump. The bed bathymetry is defined by the function \(h(x)\), with the vertical coordinate of the bed given by
\begin{equation}
h(x) = -1 + \frac{27}{4}\,\frac{H}{L^{3}}\,x\,(x-L)^{2},
\end{equation}
where we set \(H = 0.2\) and \(L = 2\) as in \cite{VanBrummelen2001}, where RANS equations are solved in the same domain for higher Reynolds numbers. The computational domain is defined as \(\Omega = [-8,\,12]\times[-1,\,0]\), with the bump positioned such that its tip lies at \(x=0\). An unstructured mesh is employed, consisting of \(46{,}731\) elements and \(801\) nodes along the free surface. Various solutions are computed using the iterative procedure described in Section \ref{sec:numericalapproach} for different combinations of Froude and Reynolds numbers, and for varying polynomial orders

\begin{figure}[ht]
    \centering
    \begin{subfigure}[b]{0.6\textwidth}
        \centering
        \includegraphics[width=\textwidth]{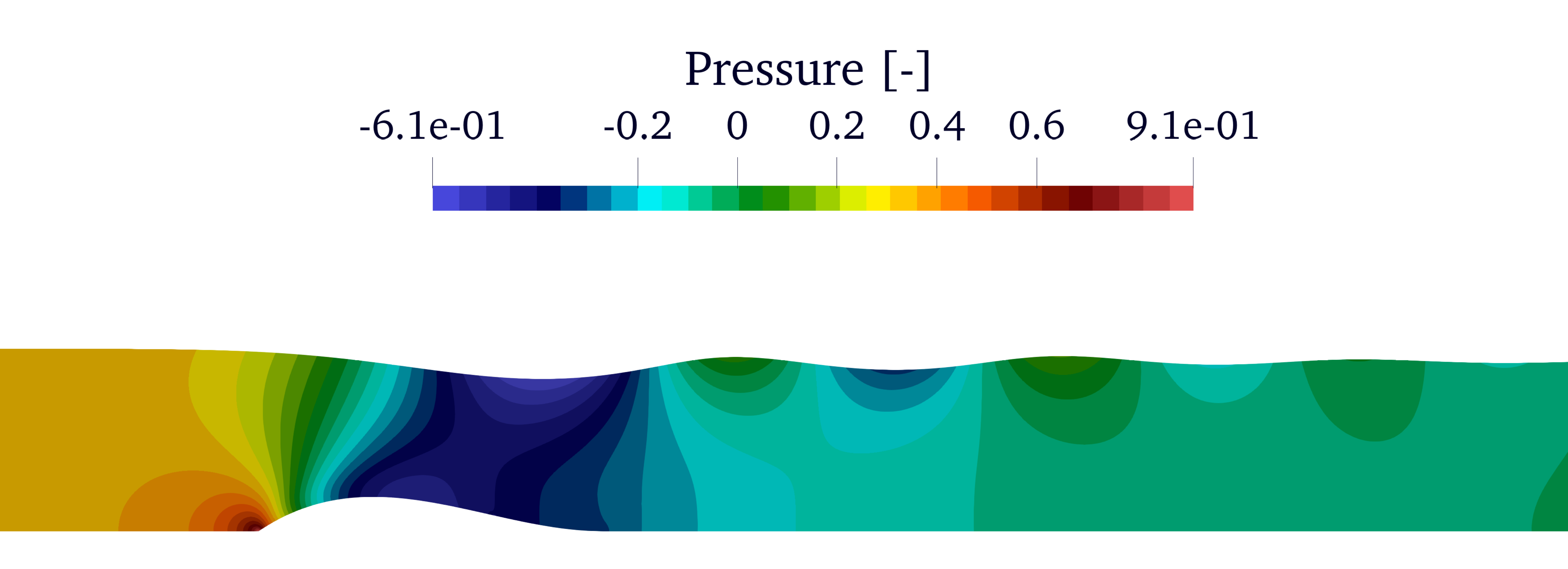}
        \caption{Dynamic Pressure.}
        \label{fig:pressure_field}
    \end{subfigure}
    \hfill
    \begin{subfigure}[b]{0.6\textwidth}
        \centering
        \includegraphics[width=\textwidth]{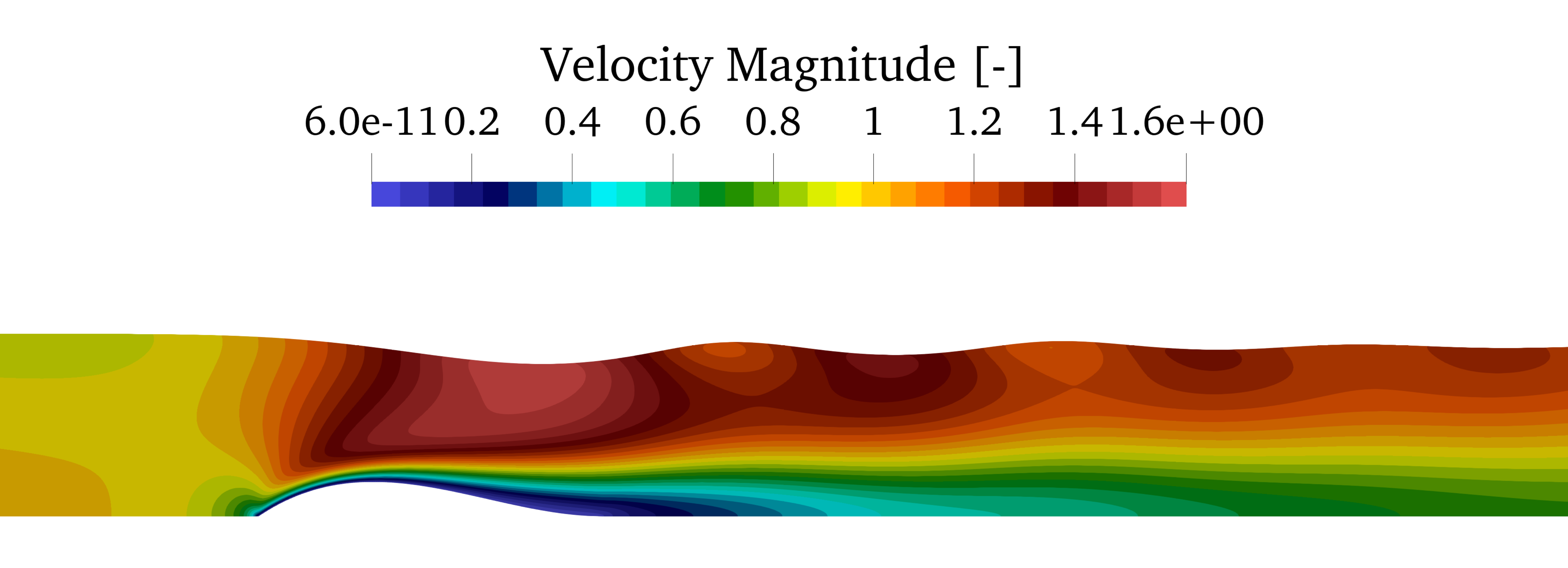}
        \caption{Velocity.}
        \label{fig:velocity_field}
    \end{subfigure}
\caption{Zoomed view of the domain showing the dynamic pressure and velocity fields for the flow over a bump with \(\mathrm{Re}=100\) and \(\mathrm{Fr}=0.43\).}
    \label{fig:pressure_velocity_re100_fr043}
\end{figure}

Figure \ref{fig:pressure_velocity_re100_fr043} reports the pressure and velocity fields for \(\mathrm{Re}=100\), \(\mathrm{Fr}=0.43\) and $P$-$Q$ = $3$-$2$. As expected, higher dynamic pressure occurs beneath wave crests, while lower values appear under troughs. In contrast, the velocity magnitude is lower below crests and higher beneath troughs. A disturbance is observed at the inflow boundary as the flow adjusts to the prescribed fixed height condition. Depending on the Froude and Reynolds numbers, the bump has a stronger or weaker influence on the upstream flow, leading to perturbations of different magnitudes.
As resolving the boundary layer developing along the bed is not the objective of this study, a free-slip boundary condition is imposed at the bottom boundary, as detailed in Section \ref{sec:mathematicalproblem}.

\subsubsection{Free surface flow validation: fully submerged NACA0012 air foil}

Next, we present the numerical results for the free surface problem involving a fully submerged NACA0012 airfoil. It is important to note that, unlike many previous studies (e.g. \cite{Lohner1999, Hino1993}), the effect of viscosity is not neglected in this work. 

\begin{figure}[ht]
    \centering
    \begin{subfigure}[b]{0.48\textwidth}
        \centering
        \includegraphics[width=\textwidth]{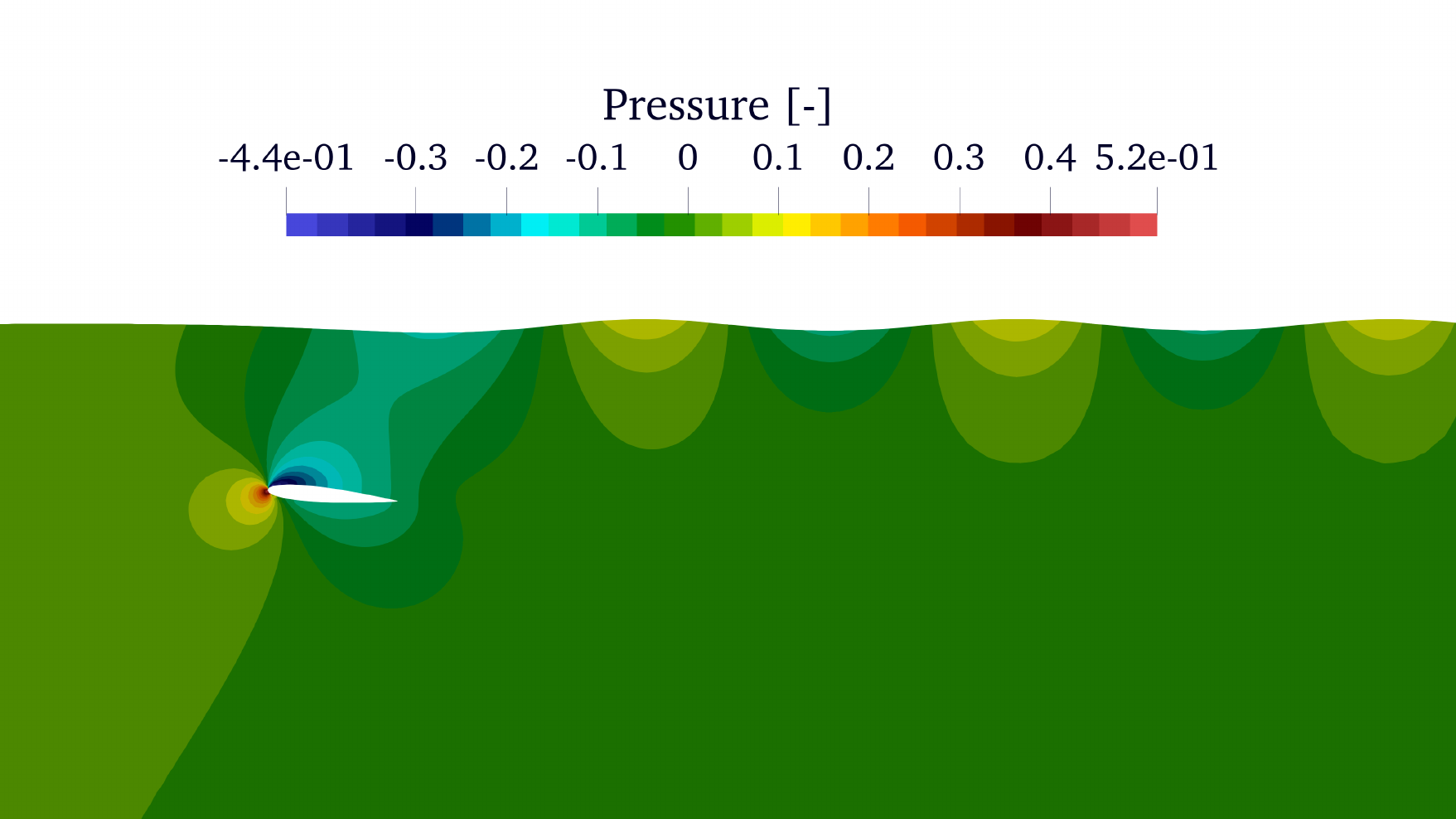}
        \caption{Dynamic pressure.}
        \label{fig:zoom_dynamic_pressure}
    \end{subfigure}
    \hfill
    \begin{subfigure}[b]{0.48\textwidth}
        \centering
        \includegraphics[width=\textwidth]{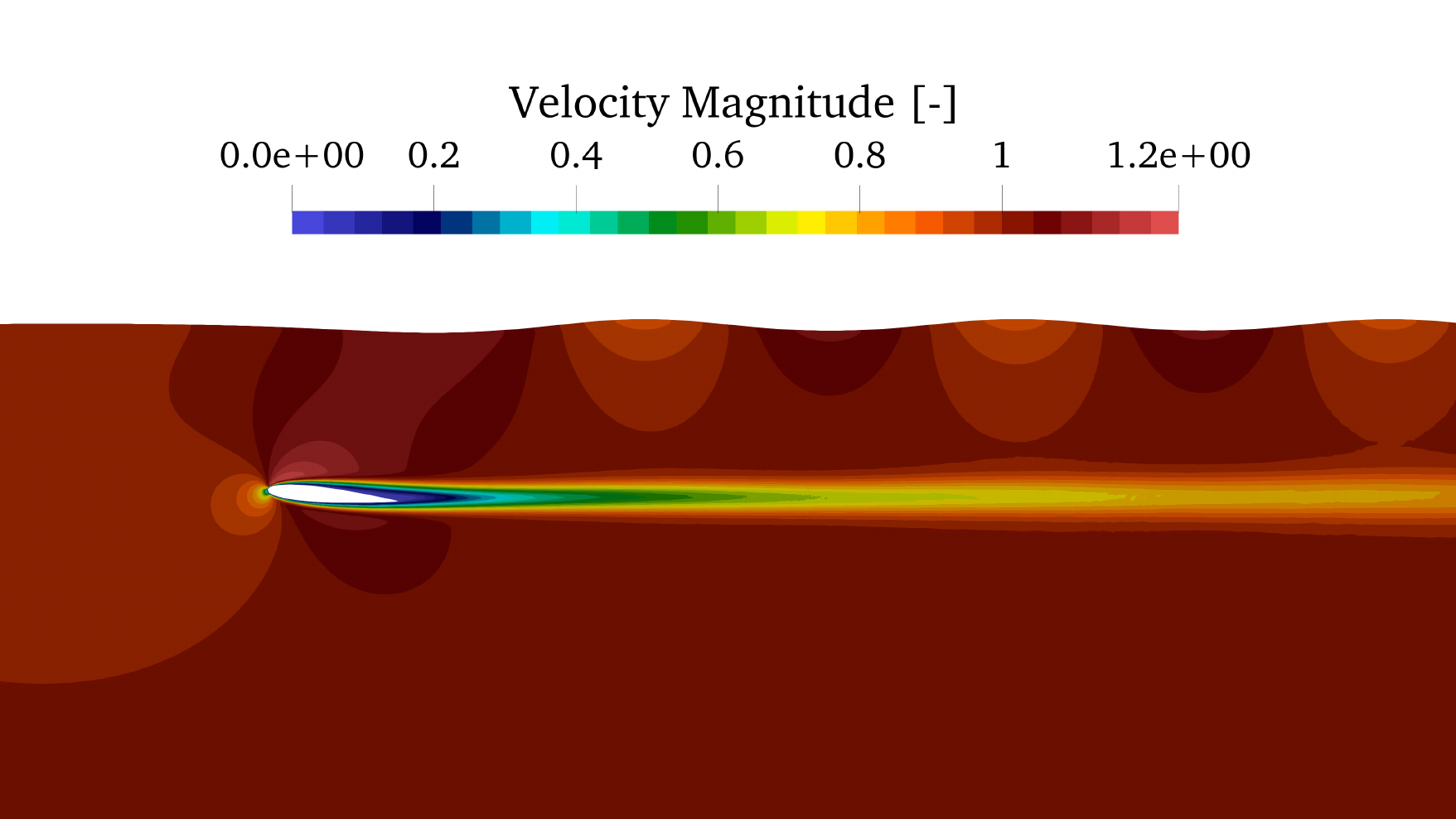}
        \caption{Velocity.}
        \label{fig:zoom_velocity}
    \end{subfigure}
\caption{Zoomed view of the domain showing the dynamic pressure and velocity fields for the free-surface flow around a NACA0012 airfoil with $\mathrm{Re}=1261$ and $\mathrm{Fr}=0.67$.}
    \label{fig:zoom_subfigures}
\end{figure}

\begin{figure}[ht]
    \centering
\includegraphics[width=0.35\textwidth]{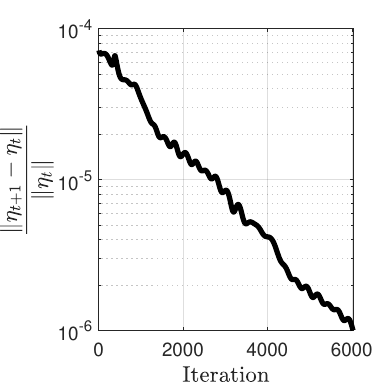}
\caption{Convergence history of the free surface elevation, showing the decay of the relative displacement as a function of the number of iterations, for the flow around a submerged NACA airfoil at \(\text{Re} = 1261\) and \(\text{Fr} = 0.67\).
    \label{fig:convergence_history}}
\end{figure}

Figure \ref{fig:zoom_dynamic_pressure} and \ref{fig:zoom_velocity} present enlarged views of a selected region in the computational domain, highlighting the detailed distributions of dynamic pressure and velocity for $\text{Re} = 1261$ and $\text{Fr} = 0.67$. The NACA airfoil is submerged at \(d = 1.252\) (scaled by the chord length), with the computational domain \(\Omega = [-14,\,20] \times [-5,\, d]\) discretized into 44{,}093 elements. The mesh resolution includes 300 edges along the free surface and 400 edges defining the airfoil boundary.
As expected, higher dynamic pressure also occurs beneath wave crests, whereas the velocity magnitude increases below troughs. Downstream of the object, the velocity field shows a wake region where the flow slows down. The plot illustrating the convergence history as a function of the number of iterations is presented in Figure \ref{fig:convergence_history}.

\begin{figure}[ht]
        \centering
        \includegraphics[width=0.8\textwidth]{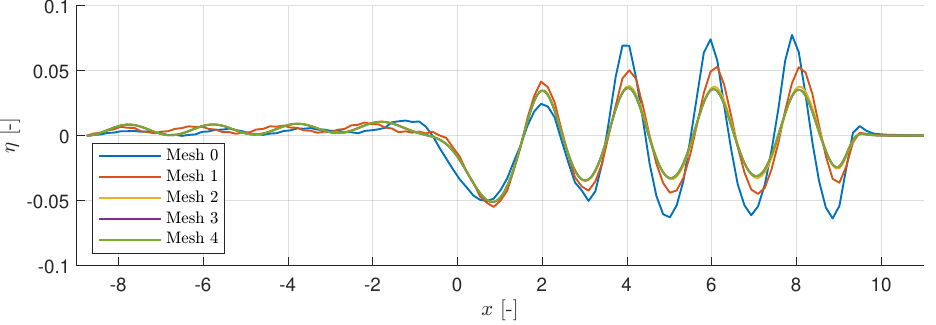}
 \caption{Surface elevation profiles for the submerged NACA0012 airfoil at $\mathrm{Re} = 1000$ and $\mathrm{Fr} = 0.57$, illustrating grid convergence across five mesh resolutions (Mesh 0: 2{,}100 elements, Mesh 1: 4{,}200 elements, Mesh 2: 15{,}500 elements, Mesh 3: 44{,}470 elements, and Mesh 4: 60{,}600 elements).}
    \label{fig:gridIndip_combined}
\end{figure}

To ensure grid independence, a study is conducted using five mesh resolutions containing from 2,100 to 60,600 elements, while keeping the polynomial degree fixed at $P$-$Q = 2$-$1$. Surface elevation profiles for $\mathrm{Re} = 1000$ are reported in Figure \ref{fig:gridIndip_combined} and show minimal variation between the finest meshes, confirming mesh convergence. Consequently, a medium-resolution mesh is chosen for the simulations to balance accuracy and computational efficiency.

\textbf{Remark.} The flow is solved for different scenarios, considering variations in Fr, Re, and depth. This visual analysis confirmed that the flow exhibits the expected characteristics under the different parameters. The results are qualitatively compared with those reported in the literature (e.g., \cite{Strom2018, duncan1983, OpenFoamDuncanPaper}); however, the corresponding graphs are omitted here, since the focus of this paper is not to document the range of possible flow regimes, but rather to present the new solver and its quantitative validation.

\textbf{Remark.} The flexibility of the solver in providing solutions for different Froude and Reynolds numbers enables the accurate reproduction of physical experiments involving steady-state fluid flows. To ensure model–prototype similarity \cite{Heller2011Scale}, the non-dimensional parameters governing the problem must be carefully selected to remain consistent with their physical counterparts. For example, to simulate a physical case corresponding to a submerged NACA profile with chord length \( c = 1 \mathrm{cm} \) and inflow velocity \( v = 0.1 \mathrm{m/s} \), moving in water with kinematic viscosity \( \nu = 10^{-6} \mathrm{m^2/s} \) under standard gravitational acceleration \( g = 9.8 \mathrm{m/s^2} \), one should consider a Reynolds number \( \mathrm{Re} = 1000 \) and a Froude number \( \mathrm{Fr} = 0.32 \). 

The influence of the inflow boundary condition, where the surface height is fixed to ensure constant fluid flux, is assessed by comparing simulations with different domain lengths at fixed $\text{Re}$ and $\text{Fr}$. Results show that extending the inflow domain has a negligible impact on the surface profile above and downstream of the submerged object, confirming the adequacy of the fixed inflow height assumption. Additionally, the effect of outflow damping is examined by comparing solutions on shorter and longer domains with consistent grid density. The damping zone effectively absorbs outgoing waves during the iterative pseudo-time step procedure. It prevents reflections at the boundary without influencing the upstream flow. Indeed, the flow solution over the object remains nearly unaffected; however, when the outflow damping is removed, non-physical disturbances arise and propagate upstream.

\subsubsection{Curvilinear verification}

\begin{figure}[H]
    \centering
    \includegraphics{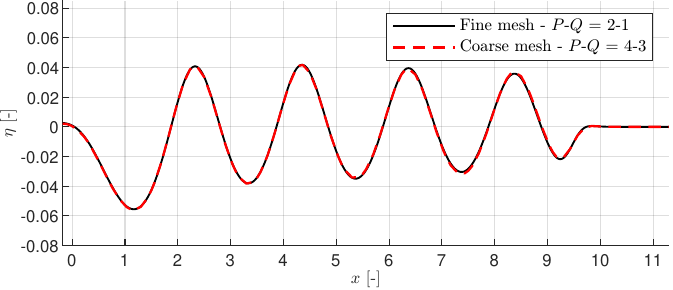}
    \caption{Comparison of mesh configurations: a fine low-order mesh with $P$-$Q$ = $2$-$1$ and 33{,}348 elements, and a coarse high-order meshes with $P$-$Q$ = $4$-$3$ and 2{,}876 elements.}
    \label{fig:curvComparison}
\end{figure}

In this section we present how curvilinear elements, combined with high-order spectral element discretization, are capable of accurately capturing the curvature of the free surface, even on relatively coarse meshes. Figure \ref{fig:curvComparison} shows a direct comparison between the results obtained using a low-order, $P$-$Q$ = $2$-$1$, fine mesh, 33{,}348 elements and 92{,}468 DOFs, and a high-order coarse mesh, with $P$-$Q$ = $4$-$3$, only 2{,}876 elements and 21{,}976 DOFs.
Despite the significant difference in mesh resolution (coarse mesh uses 50 nodes and 201 DOFs on the free surface, while the fine mesh uses 750 nodes), the two simulations exhibit nearly identical surface shapes.
This result indicates that the high-order approximation is able to accurately capture the geometry of the free surface, while using a lower number of nodes, potentially reducing the overall computational cost, as shown for the flow around a cylinder in Figure \ref{fig:ho_performance}.

\subsubsection{Comparison with \texttt{OpenFOAM} solutions: submerged NACA0012 airfoil and submerged bump}

Finally, since the physical experiments of the NACA foil and submerged bump were all performed at higher Reynolds number (e.g. \cite{duncan1983}), in this section the high-order SEM one-phase INS solver is compared to numerically obtained low Re solutions. The low Re solutions have been computed using the {\tt interFoam} transient two-phase INS solver, supplied with the well-known open-source framework {\tt OpenFOAM} (version 2212) \cite{Weller1998,openfoam} founded on the cell-centred finite volume method. The presented {\tt interFoam} solutions have all been obtained using 2nd order centred differences for gradient and diffusive terms, whereas the 2nd order TVD van Leer scheme has been used for divergence terms. The two-phase problem is solved using an algebraic volume of fluid method for the volume fraction $\alpha \in[0,\,1]$, disregarding any surface tension effects. The solutions have been time-stepped using the first-order implicit Euler scheme until approximate steady-state solutions were obtained, and the reported free surface elevations are defined by $\alpha = 0.5$.

To achieve flow similarity in the comparison of the solvers the computational domains and physical parameters are aligned, with \texttt{OpenFOAM}'s dimensional setup appropriately scaled to match the non-dimensional formulation adopted in this work. The pseudo 3D computational domains (made up of a single cell layer in the $y$-direction) consist of uniform hexahedra, with refinement regions around the free surface and the bump/foil. The baseline spatial resolution is $\Delta x = 0.003$, which was found to give "grid-independent" solutions.

\begin{figure}[H]
\centering
    \begin{subfigure}[b]{0.48\textwidth}
        \centering
        \includegraphics[width=\textwidth]{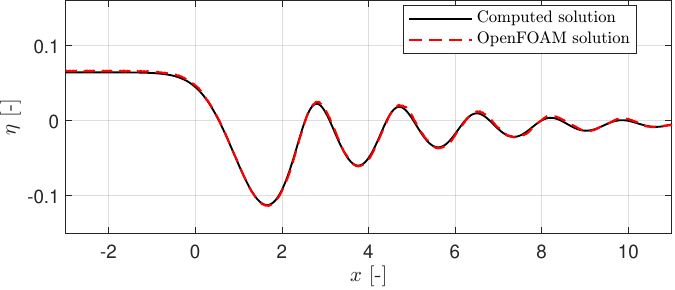}
        \caption{Flow over the submerged bump with \( \mathrm{Re} = 100 \) and \( \mathrm{Fr} = 0.43 \).}
        \label{fig:OF_bump}
    \end{subfigure}
            \hfill
    \begin{subfigure}[b]{0.48\textwidth}
        \centering
    \includegraphics[width=\textwidth]{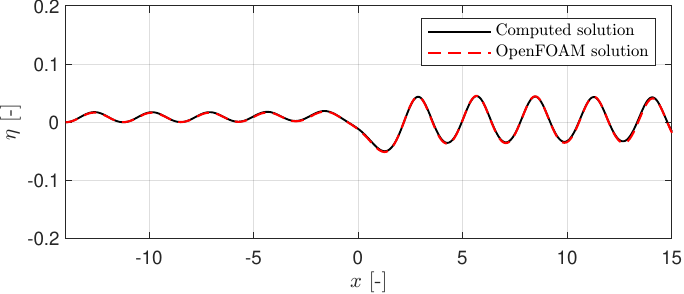}
        \caption{Flow over the submerged NACA0012 airfoil with \( \mathrm{Re} = 1261 \) and \( \mathrm{Fr} = 0.67 \).}
        \label{fig:OF_airfoil}
    \end{subfigure}
        \caption{Comparison of the steady-state free surface profiles obtained with the presented high-order solver and \texttt{OpenFOAM} for different test cases.}
        \label{fig:OFfree_surface_comparison}
\end{figure}

Figure \ref{fig:OFfree_surface_comparison} compares the steady-state free surface for the flow over the submerged bump and for the NACA airfoil case. For both cases strong agreement is observed between the results, confirming the capability of the high-order SEM method ($P$-$Q = 4$-$3$ in both comparisons) to reproduce results consistent with the reference \texttt{OpenFOAM} solution. Minor discrepancies can be attributed to the different numerical strategies and flow models employed: a single-phase formulation in the present SEM solver, in contrast to the two-phase air-fluid finite-volume approach used in \texttt{interFoam}. 

\section{Conclusion}
\label{sec:conclusion}

In this work, we present a new high-order spectral element solver for the simulation of the steady INSE subject to a free surface. This solver includes viscosity, extending previous inviscid formulations \cite{Lohner1999, Hino1993}. The mathematical and numerical model is outlined in terms of the governing equations and their implementation in \texttt{Firedrake}, together with details on mesh generation and updating (including the curvilinear element implementation) and the iterative approach used to capture the steady free surface profile.

The solver is further verified and validated through a series of numerical experiments. This includes a lid-driven cavity flow, an infinite channel flow over a fully submerged NACA0012 airfoil, and an infinite channel flow past a cylinder. All tests show good agreement with known benchmark results. In the case of the latter, spectral convergence is observed, verifying the curvilinear implementation and, moreover, justifying the use of high-order elements over low-order ones, after the accuracy and CPU time consumption are assessed. 

With the inclusion of a free surface, the solver is validated against \texttt{OpenFOAM} for finite channel flows over a bathymetry bump and a fully submerged NACA0012 airfoil. Both numerical experiments show good qualitative agreement, underlining the capability of the proposed solver to capture steady free surface profiles.

Future model development includes the natural extension to three spatial dimensions and targeting unsteady flow to capture turbulent effects.

% This work presents the development, verification, and validation of a new high-order spectral element solver for the steady INSE, with particular emphasis on free surface flows. A significant contribution of this work lies in the inclusion of viscosity, extending previous inviscid formulations as \cite{Lohner1999, Hino1993}. We first introduce the equations and the mathematical formulation of the problem, and then present the numerical implementation in \texttt{Firedrake} and its validation against classical benchmark problems. Particular focus is set on extending the solver to handle free surface flows by incorporating dynamic and kinematic boundary conditions and an iterative mesh update strategy. The methodology is tested on free surface flow around a submerged NACA airfoil, where good agreement was demonstrated compared to an \texttt{OpenFOAM} solution, underlining the solver’s ability to capture steady free surface deformation. Curvilinear elements are implemented to accurately approximate curved boundaries. 

% Future developments include extending the solver to handle partially submerged bodies in 3D, introducing unsteady capabilities, incorporating turbulence modelling, and accelerating convergence through advanced numerical strategies. These improvements would expand the solver’s capabilities, allowing it to be applied to various realistic fluid flow scenarios.

\section{Acknowledgements}

The research was carried out at the Technical University of Denmark (DTU) in the Department of Applied Mathematics and Computer Science, with computational resources provided by the DTU Computing Center (DCC) and by the National Academic Infrastructure for Supercomputing in Sweden (NAISS). Finally, the authors gratefully acknowledge the support of the Firedrake team, led by Professor David Ham. This work partly contributes to the activities of the research project led by APEK: {\em “A new digital twin concept for floating offshore structures”} supported by COWIfonden (Grant no. A-165.19) and partly to the PhD-project of JV: {\em “New Advanced Simulation Techniques for Wave Energy Converters”}.

\bibliographystyle{unsrtnat}

% When working on the draft use the following:
\bibliography{biblio.bib}  

% When submitting to arXiv (example) - remember to select all .bib etc. files.:

\end{document}